\documentclass[12pt]{article}
\usepackage{setspace}
\usepackage{amsmath,amssymb,amsthm}
\usepackage{mathrsfs}
\usepackage{geometry}
\usepackage{enumitem}
\usepackage{environ}
\usepackage{setspace}
\usepackage{titlesec}
\usepackage[all,cmtip]{xy}
\usepackage{tikz}
\usetikzlibrary{matrix,arrows}

\NewEnviron{prf}[1][]{\begin{proof}[\bf #1Proof]\BODY\end{proof}}{}
\NewEnviron{slt}[1][]{\begin{proof}[\bf #1	解]\BODY\end{proof}}{}
\newtheorem{definition}{Definition}[section]
\newtheorem{theorem}[definition]{Theorem}
\newtheorem{lemma}[definition]{Lemma}
\newtheorem{proposition}[definition]{Proposition}
\newtheorem{remark}[definition]{Remark}
\newtheorem{example}[definition]{Example}

\newtheorem{conjecture}[definition]{Conjecture}

\setlist[enumerate,1]{label=(\roman*)}
\setlist[enumerate,2]{label=(\alph*)}

\title{\textbf{On the dynamical Mordell--Lang conjecture in positive characteristic}}
\author{Junyi Xie\footnote{Beijing International Center for Mathematical Research, Beijing, China. e-mail: xiejunyi@bicmr.pku.edu.cn.}\quad and\quad She Yang\footnote{Beijing International Center for Mathematical Research, Beijing, China. e-mail: ys-yx@pku.edu.cn.}}
\date{}

\begin{document}
\begin{spacing}{1.25}

\maketitle

\begin{abstract}
We solve the dynamical Mordell--Lang conjecture for bounded-degree dynamical systems in positive characteristic. The answer in this case disproves the original version of the $p$DML conjecture.
\end{abstract}

\section{Introduction}
In this paper, as a matter of convention, every variety is assumed to be integral but the closed subvarieties can be reducible. For a rational map $f:X\dashrightarrow Y$ between two varieties, we denote $\mathrm{Dom}(f)\subseteq X$ as the domain of definition of $f$. Let $X$ be a variety over an algebraically closed field $K$ and let $f$ be a rational self-map of $X$. For a point $x\in X(K)$, we say the orbit $\mathcal{O}_{f}(x):=\{f^{n}(x)|\ n\in\mathbb{N}\}$ is well-defined if every iterate $f^{n}(x)$ lies in $\text{Dom}(f)$. We denote $\mathbb{N}=\mathbb{Z}_{+}\cup\{0\}$. An arithmetic progression is a set of the form $\{mk+l|\ k\in\mathbb{Z}\}$ for some $m,l\in\mathbb{Z}$ and an arithmetic progression in $\mathbb{N}$ is a set of the form $\{mk+l|\ k\in\mathbb{N}\}$ for some $m,l\in\mathbb{N}$.

The dynamical Mordell--Lang conjecture is one of the core problems in the field of arithmetic dynamics. It asserts that for any rational self-map $f$ of a variety $X$ over $\mathbb{C}$, the return set $\{n\in\mathbb{N}|\ f^{n}(x)\in V(\mathbb{C})\}$ is a finite union of arithmetic progressions in $\mathbb{N}$ where $x\in X(\mathbb{C})$ is a point such that the orbit $\mathcal{O}_{f}(x)$ is well-defined and $V\subseteq X$ is a closed subvariety. There is an extensive literature on various cases of this 0-DML conjecture (``0" stands for the characteristic of the base field). Two significant cases are as follows:
\begin{enumerate}
\item
If $X$ is a quasi-projective variety over $\mathbb{C}$ and $f$ is an \'etale endomorphism of $X$, then the 0-DML conjecture holds for $(X,f)$. See $\cite{Bel06}$ and $\cite[\mathrm{Theorem}\ 1.3]{BGT10}$.
\item
If $X=\mathbb{A}_{\mathbb{C}}^{2}$ and $f$ is an endomorphism of $X$, then the 0-DML conjecture holds for $(X,f)$. See $\cite{Xie17}$ and $\cite[\mathrm{Theorem}\ 3.2]{Xie}$.
\end{enumerate}

One can consult $\cite{BGT16,Xie}$ and the references therein for more known results. However, we remark that not much is known about the 0-DML conjecture when $f$ is just a \emph{rational} self-map of the variety $X$. The following problem might reflect the issue in some sense. It seems that the dynamical Mordell--Lang problem is not quite compatible with birational transformations. More precisely, let $X,Y$ be varieties over $\mathbb{C}$, $f,g$ be dominant rational self-maps of $X,Y$ respectively and $\pi:Y\dashrightarrow X$ be a birational map such that $f\circ\pi=\pi\circ g$. Even if the 0-DML conjecture holds for $(Y,g)$, generally we do not know how to deduce that the 0-DML conjecture holds for $(X,f)$.

The statement of the 0-DML conjecture fails when the base field has positive characteristic. See $\cite[\mathrm{Example}\ 3.4.5.1]{BGT16}$ for an example. Consequently, Ghioca and Scanlon proposed a dynamical Mordell--Lang conjecture in positive characteristic. See $\cite[\mathrm{Conjecture}\ 13.2.0.1]{BGT16}$. It was conjectured that the return set would be a ``$p$-normal set" (see Definition 5.1), which was firstly introduced in $\cite[\mathrm{Definition}\ 1.7]{Der07}$ for the Skolem--Mahler--Lech problem in positive characteristic. However, the Skolem--Mahler--Lech problem is a very special case of the dynamical Mordell--Lang problem (more precisely, the case of linear recurrence sequences) and it turns out that some sets of a more complicated form are needed for the general case. See Section 5 for a disproof of this old version of the $p$DML conjecture. As a result, we need to introduce some new notions of ``$p$-sets" here.

\begin{definition}
Let $p$ be a prime number and let $q=p^{e}$ for some positive integer $e$. Suppose that $d\in\mathbb{Z}_{+},r\in\mathbb{N}$ and $c_{0},c_{ij}\in\mathbb{Q}$ where $(i,j)\in\{1,\dots,d\}\times\{0,\dots,r\}$. Then we define
$$
S_{q,d,r}(c_{0};c_{10},\dots,c_{dr})=\{c_{0}+\sum\limits_{i=1}^{d}\sum\limits_{j=0}^{r} c_{ij}q^{2^{j}n_{i}}|\ n_{1},\dots,n_{d}\in\mathbb{N}\}
$$
and we shall abbreviate it as $S_{q,d,r}(c_{0};c_{ij})$ if there is no ambiguity. We define a \emph{widely} $p$\emph{-normal set in} $\mathbb{Z}$ as a union of finitely many arithmetic progressions (possibly singleton) along with finitely many \emph{subsets of} $\mathbb{Z}$ of the form $S_{q,d,r}(c_{0};c_{ij})$ as above. A \emph{widely} $p$\emph{-normal set in} $\mathbb{N}$ is a subset of $\mathbb{N}$ which is, up to a finite set, equal to the intersection of a widely $p$-normal set in $\mathbb{Z}$ and $\mathbb{N}$.
\end{definition}

Here, we say two sets $S,T$ are equal up to a finite set if the symmetric difference $(S\backslash T)\cup(T\backslash S)$ is finite, as in $\cite{Der07}$.

In the rest of this article, we will abbreviate ``widely $p$-normal set in $\mathbb{Z}$" as ``widely $p$-normal set".

\begin{remark}
\begin{enumerate}
\item
Notice that in the definition of widely $p$-normal set above, all of the sets $S_{q,d,r}(c_0;c_{ij})$ are required to be contained in $\mathbb{Z}$. We remark that this condition bounds the denominators of $c_0$ and $c_{ij}$. More precisely, we have $(q^{2^{j}}-1)\cdot\prod\limits_{s=0,s\neq j}^{r} (q^{2^{j}}-q^{2^{s}})\cdot c_{ij}\in\mathbb{Z}$ for every pair $(i,j)\in\{1,\dots,d\}\times\{0,\dots,r\}$.

\item
A widely $p$-normal set intersects $\mathbb{N}$ is a widely $p$-normal set in $\mathbb{N}$; a finite union of widely $p$-normal sets (resp. widely $p$-normal sets in $\mathbb{N}$) is still a widely $p$-normal set (resp. widely $p$-normal set in $\mathbb{N}$); $aS+b$ is a widely $p$-normal set (resp. widely $p$-normal set in $\mathbb{N}$) if $S$ is a widely $p$-normal set (resp. widely $p$-normal set in $\mathbb{N}$) and $a,b\in\mathbb{Z}$ (resp. $\mathbb{N}$). Moreover, although we will not use this statement in this article, we remark that one can show that a finite intersection of widely $p$-normal sets (resp. widely $p$-normal sets in $\mathbb{N}$) is still a widely $p$-normal set (resp. widely $p$-normal set in $\mathbb{N}$) by adapting the method of the proof of $\cite[\mathrm{Lemma}\ 9.5]{Der07}$.

\item
The definition above coincides with the definition of ``$p$-normal sets" in Definition 5.1 if one demands the $r$ always be 0.
\end{enumerate}
\end{remark}

It seems that the necessity for widely $p$-normal sets is firstly noticed in this article (see Section 5). The previous works towards the $p$DML problem are mainly as follows: $\cite{CGSZ21}$ shows that for certain cases of endomorphisms of the algebraic torus, the return set is a ``$p$-normal set in $\mathbb{N}$" (see Definition 5.1); and $\cite[\mathrm{Theorem}\ 1.4,\mathrm{Theorem}\ 1.5]{Xie23},\cite{Yang23}$ shows that the original statement of the 0-DML conjecture still holds for certain cases of high complexity endomorphisms.

People (for instance, $\cite{Xie23}$ and also Ghioca and Scanlon) guessed that the $p$-part in the return set occurs only when the dynamical system ``comes from algebraic groups". Otherwise, it is conjectured that the initial statement of the 0-DML conjecture should still be valid for $f$. We tried to form a rigorous conjecture in order to clarify this philosophy, but failed eventually. Indeed, some new phenomena and weird examples will be illustrated in an upcoming article by the authors, which show that some more complicated ``$p$-sets" will appear for certain automorphisms of algebraic groups. 

~

Now we turn to the formulation of our main result. It deals with the \emph{bounded-degree self-maps}. We once thought that $p$-sets come from bounded-degree self-maps because they are the same self-maps as those who ``come from algebraic group \emph{actions}" in some sense. But as we have mentioned above, this philosophy is not quite true.

We shall briefly introduce how to measure the complexity of a dominant rational self-map $f:X\dashrightarrow X$ in which $X$ is a projective variety over an algebraically closed field. We use the concept of the \emph{degree sequence} of $f$. Since there are many references of this concept in the literature (see for example $\cite{Dang20}$, $\cite{Tru20}$, $\cite[\mathrm{Section}\ 2.1]{Xie23}$ and $\cite[\mathrm{Section}\ 5]{Yang}$), we will just state the definition and some basic properties here.

Let $L\in\text{Pic}(X)$ be a big and nef line bundle. We consider the graph $\Gamma_f\subseteq X\times X$. Let $\pi_1,\pi_2:\Gamma_f\rightarrow X$ be the two projections. We define the \emph{first degree} $\text{deg}_{1,L}(f)$ of $f$ with respect to $L$ as the intersection number $(\pi_2^{*}(L)\cdot\pi_1^{*}(L)^{\text{dim}(X)-1})$ on $\Gamma_f$. Then we get a sequence $\{\text{deg}_{1,L}(f^n)|\ n\in\mathbb{N}\}$ of positive integers.

For two sequences $\{a_{n}\},\{b_{n}\}\in(\mathbb{R}_{\geq1})^{\mathbb{N}}$, we say that they have the same speed of growth if $\{\frac{a_{n}}{b_{n}}|\ n\in\mathbb{N}\}$ has an upper bound and a positive lower bound. Let $\mathrm{deg}_{1}(f)$ be the class of the speed of growth of the sequence $\{\mathrm{deg}_{1,L}(f^{n})|\ n\in\mathbb{N}\}$, which by $\cite[\text{Theorem}\ 1\text{(ii)}]{Dang20}$ is irrelevant with the choice of the big and nef line bundle $L$. Notice that although $\cite[\text{Theorem}\ 1\text{(ii)}]{Dang20}$ was stated for normal projective variety $X$, the result also holds for arbitrary irreducible projective variety $X$ because one can pass to the normalization of $X$. Then we can abuse notation and say that $\mathrm{deg}_{1}(f)$ is the \emph{degree sequence} of $f$. We also remark that in fact $\mathrm{deg}_{1}(f)$ remains the same on different birational models. See $\cite[\text{top of p. 1269}]{Dang20}$.

\begin{definition}
Let $X$ be a projective variety over an algebraically closed field. We say a dominant rational self-map $f:X\dashrightarrow X$ is of bounded-degree if $\mathrm{deg}_{1}(f)$ is a bounded sequence.
\end{definition}

Notice that according to our definition, a bounded-degree self-map is a priori dominant.

\begin{remark}
If $f$ is a surjective endomorphism of a projective variety $X$ in the definition above, then $f$ is of bounded-degree if and only if the sequence $\{(f^{n})^{*}(L)\cdot L^{\mathrm{dim}(X)-1}|\ n\in\mathbb{N}\}$ is bounded for some (and hence every) ample line bundle $L\in\mathrm{Pic}(X)$.
\end{remark}

Now we can state the main theorem of this paper.

\begin{theorem}
Let $X$ be a projective variety over an algebraically closed field $K$ of characteristic $p>0$ and let $f$ be a bounded-degree self-map of $X$. Let $x\in X(K)$ be a closed point such that the orbit $\mathcal{O}_{f}(x)$ is well-defined and let $V\subseteq X$ be a closed subvariety. Then $\{n\in\mathbb{N}|\ f^{n}(x)\in V(K)\}$ is a widely $p$-normal set in $\mathbb{N}$.
\end{theorem}

\begin{remark}
\begin{enumerate}
\item
We only consider bounded-degree self-maps of projective varieties in this article, but one can also define this notion on quasi-projective varieties since the degree sequence is independent of the choice of birational models as we have mentioned before. Then our result automatically extends to bounded-degree self-maps of quasi-projective varieties over an algebraically closed field of positive characteristic.

\item
We focus on the positive characteristic case in this article, but our method is also valid for the 0-characteristic case. Notice $\cite[\mathrm{Theorem}\ 7]{CS93}$ implies that $\{n\in\mathbb{Z}|\ g^{n}\in X(K)\}$ is a finite union of arithmetic progressions in which $G$ is an algebraic group over an algebraically closed field $K$ of characteristic 0, $g\in G(K)$ is a closed point, and $X\subseteq G$ is a closed subvariety. By using this result instead of Theorem 3.1 in the proof of Theorem 1.5, one can show that the set $\{n\in\mathbb{N}|\ f^{n}(x)\in V(K)\}$ is a finite union of arithmetic progressions in $\mathbb{N}$ if $K$ is an algebraically closed field of characteristic 0 (and the other conditions remain the same).
\end{enumerate}
\end{remark}

At the end of the Introduction, we describe the structure of this paper. In Section 2, we make some preparations about the Mordell--Lang problem in positive characteristic. Then, in Section 3, we use those preparations to prove Theorem 1.5 in the case of translation of algebraic groups. After that, we finish the proof of Theorem 1.5 in Section 4, using the philosophy that bounded-degree self-maps come from group actions. More concretely, there is a regularization theorem (Theorem 4.6) which says that a bounded-degree self-map can be regularized into a bounded-degree automorphism of a projective variety after some birational transformation. Then by noticing that such bounded-degree automorphisms are precisely those automorphisms which come from algebraic group actions, we may prove Theorem 1.5 by using Theorem 3.1. Lastly, in Section 5, we give some heuristic examples in which the widely $p$-normal sets may appear in the return set and provide a rigorous disproof of the original version of the $p$DML conjecture.

This paper also includes an Appendix about the $p$-Mordell--Lang problem in which we fix an error and generalize the results in $\cite{GY}$.

\section{Preparations}
In this section, we will make some preparations about the Mordell--Lang problem in positive characteristic. These preparations will be used in Section 3 to prove that the statement of Theorem 1.5 holds for translations of algebraic groups.

\subsection{Two theorems towards Mordell--Lang problem in positive characteristic}
In this subsection, we introduce the main theorems of $\cite{Hru96}$ and $\cite{MS04}$ towards the Mordell--Lang problem in positive characteristic.

The Mordell--Lang conjecture for a semiabelian variety $G$ over an algebraically closed field $K$ of characteristic 0 states that $X(K)\cap\Gamma$ is a finite union of cosets of subgroups of $\Gamma$ where $X\subseteq G$ is a closed subvariety and $\Gamma\subseteq G(K)$ is a finitely generated subgroup. This conjecture was proved by Faltings $\cite{Fal94}$ in the case of abelian varieties and by Vojta $\cite{Voj96}$ in the general case of semiabelian varieties. But this statement fails in positive characteristic.

In the positive characteristic case, the groundbreaking work $\cite[\mathrm{Theorem}\ 1.1]{Hru96}$ proves that the counterexamples of the original 0-ML statement in positive characteristic come from isotrivial semiabelian varieties. Then $\cite[\mathrm{Theorem}\ \mathrm{B}]{MS04}$ describes the form of the intersection set of a closed subvariety and an Frobenius invariant finitely generated subgroup in an isotrivial semiabelian variety. Recently, $\cite{Ghi}$ gives a description of the intersection set of a closed subvariety and an arbitrary finitely generated subgroup in an isotrivial semiabelian variety.

We firstly state the main theorem of $\cite{Hru96}$.

\begin{theorem}($\cite[\mathrm{Theorem}\ 1.1]{Hru96}$)
Let $G$ be a semiabelian variety over an algebraically closed field $K$ of characteristic $p>0$. Let $X\subseteq G$ be an irreducible closed subvariety and let $\Gamma\subseteq G(K)$ be a finitely generated subgroup. If $X(K)\cap\Gamma$ is dense in $X$, then there exist

$\bullet$ a semiabelian subvariety $G_{1}\subseteq G$,

$\bullet$ a semiabelian variety $G_{0}$ over a finite subfield $\mathbb{F}_{q}\subseteq K$,

$\bullet$ a geometrically integral closed subvariety $X_{0}\subseteq G_{0}$,

$\bullet$ a surjective algebraic group homomorphism $f:G_{1}\rightarrow G_{0}\times_{\mathbb{F}_{q}}K$, and

$\bullet$ a point $x_{0}\in G(K)$

such that $X=x_{0}+f^{-1}(X_{0}\times_{\mathbb{F}_{q}}K)$ as a closed subset of $G$.
\end{theorem}

Now we state $\cite[\mathrm{Theorem}\ \mathrm{B}]{MS04}$ concerning the form of the intersection set of a closed subvariety and an Frobenius invariant finitely generated subgroup in an isotrivial semiabelian variety. Firstly, we introduce the definition of ``$F$-sets".

\begin{definition}
Let $G$ be a semiabelian variety over an algebraically closed field $K$ of characteristic $p>0$ which is defined over a finite field $\mathbb{F}_{q}$. Let $F=\mathrm{Frob}_{q}$ be the Frobenius endomorphism of $G$ and let $\Gamma\subseteq G(K)$ be a finitely generated subgroup.
\begin{enumerate}
\item
An $F$-\emph{set} is a subset of $G(K)$ of the form $\{\alpha_0+\sum\limits_{i=1}^{d}F^{kn_i}(\alpha_{i})|\ n_1,\dots,n_d\in\mathbb{N}\}$ where $k,d\in\mathbb{Z}_{+}$ and $\alpha_{0},\alpha_{1},\dots,\alpha_{d}\in G(K)$. An $F$-\emph{set in} $\Gamma$ is an $F$-set which is contained in $\Gamma$.

\item
An $F$-\emph{normal set in} $\Gamma$ is a finite union of sets of the form $S+\Lambda$ where $S$ is an $F$-set in $\Gamma$ and $\Lambda\subseteq\Gamma$ is a subgroup.
\end{enumerate}
\end{definition}

\begin{remark}
In literature, the ``$F$-set" above may be called as ``groupless $F$-set" and the ``$F$-normal set" above may be called as ``$F$-set". We use this slightly modified terminology here in order to be similar as in Definition 1.1. Moreover, we would like to shorten the name of those $F$-sets in Definition 2.2(i) because we will mainly play with them later.
\end{remark}

\begin{theorem}($\cite[\mathrm{Theorem}\ \mathrm{B}]{MS04}$)
Let $G$ be a semiabelian variety over an algebraically closed field $K$ of characteristic $p>0$ which is defined over a finite field $\mathbb{F}_{q}$. Let $F=\mathrm{Frob}_{q}$ be the Frobenius endomorphism of $G$. Let $\Gamma\subseteq G(K)$ be a finitely generated $F$-\emph{invariant} subgroup and let $X\subseteq G$ be a closed subvariety. Then $X(K)\cap\Gamma$ is an $F$-normal set in $\Gamma$.
\end{theorem}

Notice that the $F$-invariant requirement on the finitely generated subgroup $\Gamma$ (i.e. $F(\Gamma)\subseteq\Gamma$) in the theorem above cannot be dropped.

\subsection{Intersection of $F$-sets and cyclic subgroups}
As we have mentioned above, Theorem 2.4 is only available for the intersection of closed subvarieties with $F$-invariant finitely generated subgroups. But for our application towards the $p$DML problem, we mainly concern about the intersection set of a closed subvariety with a cyclic group. So we will discuss the form of the intersection set of an $F$-set with an infinite cyclic group in this subsection.

We start with the definition of ``widely $F$-sets". Their form is similar to the form of ``widely $p$-sets" in Definition 1.1.

\begin{definition}
Let $G$ be a semiabelian variety over an algebraically closed field $K$ of characteristic $p>0$ which is defined over a finite field $\mathbb{F}_{q}$. Let $F=\mathrm{Frob}_{q}$ be the Frobenius endomorphism of $G$. A \emph{widely} $F$-\emph{set} is a subset of $G(K)$ of the form $\{\alpha_0+\sum\limits_{i=1}^{d}\sum\limits_{j=0}^{r} F^{k\cdot2^{j}n_i}(\alpha_{ij})|\ n_1,\dots,n_d\in\mathbb{N}\}$ where $k,d\in\mathbb{Z}_{+},r\in\mathbb{N}$ and $\alpha_{0},\alpha_{ij}\in G(K)$ where $(i,j)\in\{1,\dots,d\}\times\{0,\dots,r\}$.
\end{definition}

The main result of this subsection is the following proposition.

\begin{proposition}
Let $G$ be a semiabelian variety over an algebraically closed field $K$ of characteristic $p>0$ which is defined over a finite field $\mathbb{F}_{q}$. Let $F=\mathrm{Frob}_{q}$ be the Frobenius endomorphism of $G$ and let $\Gamma\subseteq G(K)$ be a cyclic group. Then the intersection of an $F$-set with $\Gamma$ is a finite union of widely $F$-sets.
\end{proposition}

\begin{remark}
\begin{enumerate}
\item
We will see that the common ratio 2 of the geometric series on the exponents comes from the two different absolute values of the roots of the minimal polynomial of the Frobenius endomorphism.

\item
With some efforts, one can show that the conclusion above also holds for the intersection of a (widely) $F$-set with a finitely generated subgroup of $G(K)$. But since we only need this general statement in the Appendix and its proof is essentially the same as the proof of Proposition 2.6, we will focus on the intersection of an $F$-set with a cyclic group for simplicity.
\end{enumerate}
\end{remark}

In order to make things clear, we start with a technical definition.

\begin{definition}
Let $d\in\mathbb{Z}_+$.
\begin{enumerate}
\item
A \emph{good subgroup of} $\mathbb{Z}^{d}=\{(n_1,\dots,n_d)|\ n_{1},\dots,n_{d}\in\mathbb{Z}\}$ is a subgroup defined by some requirements of the following 4 types
\begin{enumerate}
\item
$n_{i}=0$,
\item
$n_{i}$ is a multiple of some positive integer $D$,
\item
$n_{i}=n_{j}$, and
\item
$n_{i}=2\cdot n_{j}$
\end{enumerate}
where $1\leq i,j\leq d$.

\item
A \emph{big rectangle in} $\mathbb{N}^{d}$ is a set of the form $\{(n_1,\dots,n_d)\in\mathbb{N}^{d}|\ n_{1}\geq M_{1},\dots,n_{d}\geq M_{d}\}$ for some nonnegative integers $M_{1},\dots,M_{d}$.

\item
A \emph{good coset in} $\mathbb{N}^{d}$ is the intersection set of a big rectangle in $\mathbb{N}^{d}$ with a coset of a good subgroup of $\mathbb{Z}^{d}$.
\end{enumerate}
\end{definition}

The following lemma is easy to verify and we omit its proof.

\begin{lemma}
\begin{enumerate}
\item
The intersection of two good subgroups of $\mathbb{Z}^{d}$ is also a good subgroup of $\mathbb{Z}^{d}$; the intersection of two good cosets in $\mathbb{N}^{d}$ is also a good coset in $\mathbb{N}^{d}$.

\item
A good subgroup of $\mathbb{Z}^{d}$ can be expressed as the form $\sum\limits_{i=1}^{k}\mathbb{Z}\cdot D_{i}\eta_{i}$ where $k\in\mathbb{N}$, $D_1,\dots,D_{k}\in\mathbb{Z}_{+}$ and $\eta_1,\dots,\eta_{k}\in\mathbb{N}^{d}$ are vectors satisfying
\begin{enumerate}
\item
every component of each $\eta_{i}$ lies in $\{0\}\cup\{2^{m}|\ m\in\mathbb{N}\}$, and

\item
the vectors $\eta_{i}$ are pairwise orthogonal, i.e. any two of them do not have a common nonzero component.
\end{enumerate}
A nonempty good coset in $\mathbb{N}^{d}$ can be expressed as the form $\eta_{0}+\sum\limits_{i=1}^{k}\mathbb{N}\cdot D_{i}\eta_{i}$ where $\eta_{0}\in\mathbb{N}^{d}$ and $k,D_1,\dots,D_{k},\eta_1,\dots,\eta_{k}$ are as above.
\end{enumerate}
\end{lemma}

\begin{remark}
In fact, the vectors $\eta_{i}$ in part (ii) above also have a property that their components are ``continuous", i.e. $\eta_{i}$ will have a component which is $2^{a}$ if it has a component which is $2^{b}$ for some $b>a$.
\end{remark}

Now we can prove Proposition 2.6.

\proof[Proof of Proposition 2.6]
Let $\Gamma\subseteq G(K)$ be a cyclic subgroup generated by $g\in G(K)$. We may assume that $g$ is non-torsion because otherwise $\Gamma$ will be a finite group and hence the result follows. Recall that an $F$-set is a set of the form $\{\alpha_0+\sum\limits_{i=1}^{d}F^{kn_i}(\alpha_{i})|\ n_1,\dots,n_d\in\mathbb{N}\}$ where $k,d\in\mathbb{Z}_{+}$ and $\alpha_{0},\alpha_{1},\dots,\alpha_{d}\in G(K)$. We denote $\Phi=F^{k}$ and thus $\Phi$ is also a Frobenius endomorphism of $G$, namely, $\Phi=\mathrm{Frob}_{q^{k}}$. In the following, we will prove that if the intersection set of $\{\alpha_0+\sum\limits_{i=1}^{d}\Phi^{n_i}(\alpha_{i})|\ n_1,\dots,n_d\in\mathbb{N}\}$ and $\Gamma$ is infinite, then $\{(n_1,\dots,n_{d})\in\mathbb{N}^{d}|\ \alpha_0+\sum\limits_{i=1}^{d}\Phi^{n_i}(\alpha_{i})\in\Gamma\}$ is a finite union of good cosets in $\mathbb{N}^{d}$ and hence the result follows. At this point, we would like to make two remarks: firstly, one can let some $\alpha_{ij}$ be $0$ in the definition of widely $F$-sets; secondly, although the positive integers $D_i$ in Lemma 2.9(ii) can be different, one may consider their least common multiple and then split a set of the form $\{\alpha_0+\sum\limits_{i=1}^{d}\sum\limits_{j=0}^{r} F^{kD_{i}\cdot2^{j}n_i}(\alpha_{ij})|\ n_1,\dots,n_d\in\mathbb{N}\}$ into finitely many widely $F$-sets.

Now firstly, recall that $\Phi$ admits an equation $P(\Phi)=0$ where $P(x)=(x-a_{1})\cdots(x-a_{s})\in\mathbb{Z}[x]$ is a monic $\mathbb{Z}$-coefficient polynomial whose complex roots $a_{1},\dots,a_{s}$ are pairwise distinct and all have absolute value $q^{k}$ or $q^{\frac{k}{2}}$. For each $0\leq j\leq s-1$, we let $c_{j}(n)$ be the linear recurrence sequence satisfying $c_{j}(j)=1$, $c_{j}(0)=\cdots=c_{j}(j-1)=c_{j}(j+1)\cdots=c_{j}(s-1)=0$, and has $P(x)$ as its characteristic polynomial. Then for each $0\leq j\leq s-1$, there exist $b_{j1},\dots,b_{js}\in\mathbb{C}$ such that $c_{j}(n)=b_{j1}a_{1}^{n}+\cdots+b_{js}a_{s}^{n}$ for every $n\in\mathbb{N}$. Moreover, we have $\Phi^{n}(\alpha)=\sum\limits_{j=0}^{s-1} c_{j}(n)\Phi^{j}(\alpha)$ for every nonnegative integer $n$.

Let $M$ be the finitely generated subgroup of $G(K)$ generated by $\{\alpha_0\}\cup\{\Phi^{j}(\alpha_{i})|\ 1\leq i\leq d,0\leq j\leq s-1\}$. So $\{\alpha_0+\sum\limits_{i=1}^{d}\Phi^{n_i}(\alpha_{i})|\ n_1,\dots,n_d\in\mathbb{N}\}\subseteq M$. If $M\cap\Gamma=\{0\}$, then the intersection set of $\{\alpha_0+\sum\limits_{i=1}^{d}\Phi^{n_i}(\alpha_{i})|\ n_1,\dots,n_d\in\mathbb{N}\}$ and $\Gamma$ is contained in $\{0\}$ and hence is a finite union of widely $F$-sets. So we may assume that $M\cap\Gamma$ is an infinite cyclic group generated by $g_{0}\in M$. We want to show that $\{(n_1,\dots,n_{d})\in\mathbb{N}^{d}|\ \alpha_0+\sum\limits_{i=1}^{d}\Phi^{n_i}(\alpha_{i})\in\mathbb{Z}\cdot g_{0}\}$ is a finite union of good cosets in $\mathbb{N}^{d}$.

Write $M=M_{\mathrm{tor}}\oplus M_{0}$ where $M_{0}$ is a free $\mathbb{Z}$-module. We can find positive integers $m_{0},D$ and a $\mathbb{Z}$-basis $\{R_{1},\dots,R_{l}\}$ of $M_{0}$ such that $m_{0}\cdot g_{0}=D\cdot R_{1}$. Now fix any integer $r\in\{0,\dots,m_{0}-1\}$, we will prove that $\{(n_1,\dots,n_{d})\in\mathbb{N}^{d}|\ \alpha_0+\sum\limits_{i=1}^{d}\Phi^{n_i}(\alpha_{i})\in(r+m_{0}\mathbb{Z})\cdot g_{0}\}$ is a finite union of good cosets in $\mathbb{N}^{d}$.

Write $\alpha_{0}-r\cdot g_{0}=T_{0}+\sum\limits_{u=1}^{l}x_{u}\cdot R_{u}$ and $\Phi^{j}(\alpha_{i})=T_{ij}+\sum\limits_{u=1}^{l}y_{iju}\cdot R_{u}$ for some $T_{0},T_{ij}\in M_{\mathrm{tor}}$ and $x_{u},y_{iju}\in\mathbb{Z}$ where $1\leq i\leq d,0\leq j\leq s-1$ and $1\leq u\leq l$. Then for each $1\leq i\leq d$, we have $\Phi^{n_i}(\alpha_{i})=\sum\limits_{j=0}^{s-1}c_{j}(n_i)\Phi^{j}(\alpha_i)=\sum\limits_{j=0}^{s-1}c_{j}(n_i)(T_{ij}+\sum\limits_{u=1}^{l}y_{iju}\cdot R_{u})$. So\\
\begin{equation*}
\begin{split}
(\alpha_0-r\cdot g_{0})+\sum\limits_{i=1}^{d}\Phi^{n_i}(\alpha_{i})=(T_{0}+\sum\limits_{u=1}^{l}x_{u}\cdot R_{u})+\sum\limits_{i=1}^{d}\sum\limits_{j=0}^{s-1}c_{j}(n_i)(T_{ij}+\sum\limits_{u=1}^{l}y_{iju}\cdot R_{u})\\=(T_0+\sum\limits_{i=1}^{d}\sum\limits_{j=0}^{s-1}c_{j}(n_i)\cdot T_{ij})+\sum\limits_{u=1}^{l}(x_u+\sum\limits_{i=1}^{d}\sum\limits_{j=0}^{s-1}c_{j}(n_i)y_{iju})\cdot R_u.
\end{split}
\end{equation*}
Therefore, we may interpret the condition on $(n_1,\dots,n_{d})\in\mathbb{N}^{d}$ that $\alpha_0+\sum\limits_{i=1}^{d}\Phi^{n_i}(\alpha_{i})\in(r+m_{0}\mathbb{Z})\cdot g_{0}$ as follows (recall that $m_{0}\cdot g_{0}=D\cdot R_{1}$):
\begin{enumerate}
\item
$T_0+\sum\limits_{i=1}^{d}\sum\limits_{j=0}^{s-1}c_{j}(n_i)\cdot T_{ij}=0$,

\item
$x_1+\sum\limits_{i=1}^{d}\sum\limits_{j=0}^{s-1}c_{j}(n_i)y_{ij1}$ is a multiple of $D$, and

\item
$x_u+\sum\limits_{i=1}^{d}\sum\limits_{j=0}^{s-1}c_{j}(n_i)y_{iju}=0$ for every $2\leq u\leq l$.
\end{enumerate}

Now we argue that for each one of those three conditions, the set $\{(n_1,\dots,n_d)\in\mathbb{N}^{d}|\ (n_1,\dots,n_d)$ satisfies that condition$\}$ is a finite union of good cosets in $\mathbb{N}^{d}$. Then by Lemma 2.9(i), we conclude that the intersection set $\{(n_1,\dots,n_d)\in\mathbb{N}^{d}|\ \alpha_0+\sum\limits_{i=1}^{d}\Phi^{n_i}(\alpha_{i})\in(r+m_{0}\mathbb{Z})\cdot g_{0}\}$ is also a finite union of good cosets in $\mathbb{N}^{d}$.

Firstly, we deal with condition (i). Notice that for every $0\leq j\leq s-1$, the sequence $(c_{j}(n)\ (\mathrm{mod}\ |M_{\mathrm{tor}}|))_{n\in\mathbb{N}}$ is eventually periodic. So for every $1\leq i\leq d$, the sequence of torsion points $(\sum\limits_{j=0}^{s-1}c_{j}(n)\cdot T_{ij})_{n\in\mathbb{N}}$ is eventually periodic. Now one can easily see that the set $\{(n_1,\dots,n_d)\in\mathbb{N}^{d}|\ (n_1,\dots,n_d)$ satisfies condition (i)$\}$ is a finite union of good cosets in $\mathbb{N}^{d}$. Indeed, the good subgroups of $\mathbb{Z}^{d}$ involved here can be defined using only requirements of type (a)(b) in Definition 2.8(i).

Secondly, for condition (ii), the argument is just the same as in that for condition (i). One just substitute ``$\mathrm{mod}\ |M_{\mathrm{tor}}|$" by ``$\mathrm{mod}\ D$".

Lastly, we deal with condition (iii). Using Lemma 2.9(i), we only need to prove that the set $\{(n_1,\dots,n_d)\in\mathbb{N}^{d}|\ x_u+\sum\limits_{i=1}^{d}\sum\limits_{j=0}^{s-1}c_{j}(n_i)y_{iju}=0\}$ is a finite union of good cosets in $\mathbb{N}^{d}$ for every fixed integer $u\in\{2,\dots,l\}$. So now we fix $u$ and calculate that\\
$$
x_u+\sum\limits_{i=1}^{d}\sum\limits_{j=0}^{s-1}c_{j}(n_i)y_{iju}=x_{u}+\sum\limits_{i=1}^{d}\sum\limits_{j=0}^{s-1}y_{iju}(\sum\limits_{v=1}^{s}b_{jv}a_{v}^{n_i})=x_u+\sum\limits_{i=1}^{d}\sum\limits_{v=1}^{s}(\sum\limits_{j=0}^{s-1}y_{iju}b_{jv})a_{v}^{n_i}.
$$
To be clear, we let $C_{iv}=\sum\limits_{j=0}^{s-1}y_{iju}b_{jv}$ and then the condition can be written as $x_u+\sum\limits_{i=1}^{d}\sum\limits_{v=1}^{s}C_{iv}a_{v}^{n_i}=0$ in which $C_{iv}\in\mathbb{C}$ are constants for $(i,v)\in\{1,\dots,d\}\times\{1,\dots,s\}$.

Now consider the elements $\lambda_{i}=(1,\dots,1,a_{1},\dots,a_{s},1,\dots,1)\in\mathbb{G}_{m}^{ds}(\mathbb{C})$ for $1\leq i\leq d$ in which $a_{s}$ lies at the $is$th coordinate in the vector $\lambda_i$. Let $V\subseteq\mathbb{G}_{m,\mathbb{C}}^{ds}$ be the hyperplane defined by the linear equation $x_{u}+\sum\limits_{i=1}^{d}\sum\limits_{v=1}^{s}C_{iv}X_{(i-1)s+v}=0$ where $X_1,\dots,X_{ds}$ are the coordinates of $\mathbb{G}_{m,\mathbb{C}}^{ds}$. Then for $(n_1,\dots,n_d)\in\mathbb{N}^{d}$, we interpret the condition $x_u+\sum\limits_{i=1}^{d}\sum\limits_{j=0}^{s-1}c_{j}(n_i)y_{iju}=0$ as $\sum\limits_{i=1}^{d}n_{i}\cdot\lambda_{i}\in V(\mathbb{C})$. Let $\Lambda\subseteq\mathbb{G}_{m}^{ds}(\mathbb{C})$ be the free subgroup with $\mathbb{Z}$-basis $\{\lambda_1,\dots,\lambda_d\}$. Then $\cite[\text{Th\'eor\`eme}\ 2]{Lau84}$ says that $\{\sum\limits_{i=1}^{d}n_{i}\cdot\lambda_{i}|\ n_{1},\dots,n_{d}\in\mathbb{Z}\ \text{and}\ \sum\limits_{i=1}^{d}n_{i}\cdot\lambda_{i}\in V(\mathbb{C})\}$ is a finite union of sets of the form $\lambda+(H\cap\Lambda)$ in which $\lambda\in\Lambda$ and $H\subseteq\mathbb{G}_{m,\mathbb{C}}^{ds}$ is an algebraic subgroup defined by some equations of the form $X_{i}=1$ or $X_{i}=X_{j}$ where $1\leq i,j\leq ds$. Now recall that $a_1,\dots,a_s$ are complex numbers having absolute value $q^{k}$ or $q^{\frac{k}{2}}$, we conclude that $\{(n_1,\dots,n_d)\in\mathbb{Z}^{d}|\ \sum\limits_{i=1}^{d}n_{i}\cdot\lambda_{i}\in V(\mathbb{C})\}$ is a finite union of cosets of good subgroups of $\mathbb{Z}^{d}$. As a result, we have proved that $\{(n_1,\dots,n_d)\in\mathbb{N}^{d}|\ x_u+\sum\limits_{i=1}^{d}\sum\limits_{j=0}^{s-1}c_{j}(n_i)y_{iju}=0\}$ is a finite union of good cosets in $\mathbb{N}^{d}$ for every fixed integer $u\in\{2,\dots,l\}$.

To sum up, the whole argument above proves that $\{(n_1,\dots,n_{d})\in\mathbb{N}^{d}|\ \alpha_0+\sum\limits_{i=1}^{d}\Phi^{n_i}(\alpha_{i})\in\Gamma\}$ is a finite union of good cosets in $\mathbb{N}^{d}$. So taking the argument at the end of the first paragraph in the proof into account, we finish the proof of Proposition 2.6.
\endproof

\subsection{The form of widely $F$-sets contained in an infinite cyclic group}
In subsection 2.2, we have proved that the intersection of an $F$-set with a cyclic group is a finite union of widely $F$-sets. In this subsection, we will discuss the form of those widely $F$-sets contained in a cyclic group.

Let $G$ be a semiabelian variety over an algebraically closed field $K$ of characteristic $p>0$ which is defined over a finite field $\mathbb{F}_{q}$. Let $F=\mathrm{Frob}_{q}$ be the Frobenius endomorphism of $G$ and let $\Gamma\subseteq G(K)$ be an infinite cyclic group. Since our purpose is to describe the form of the intersection of $\Gamma$ with a closed subvariety $X\subseteq G$, it is natural to consider that whether $\Gamma$ has a nontrivial $F^{k}$-invariant subgroup for some $k\in\mathbb{Z}_{+}$. If so, then we may apply Theorem 2.4 and hence in this case, we only need to consider the form of $F$-sets contained in an infinite cyclic group. The result to this end is as follows.

\begin{proposition}
Let $G$ be a semiabelian variety over an algebraically closed field $K$ of characteristic $p>0$ which is defined over a finite field $\mathbb{F}_{q}$. Let $F=\mathrm{Frob}_{q}$ be the Frobenius endomorphism of $G$. Let $\Gamma$ be the cyclic group generated by a non-torsion point $g\in G(K)$ and let $S$ be an $F$-set contained in $\Gamma$. Then the set $\{n\in\mathbb{Z}|\ n\cdot g\in S\}$ is a finite union of sets of the form $S_{p^{m},d,0}(c_{0};c_{i0})$ defined in Definition 1.1. In other words, the ``widely" $p$-sets will not appear in here.
\end{proposition}

\begin{prf}
Firstly, recall that $S$ is of the form $\{\alpha_0+\sum\limits_{i=1}^{d} F^{kn_{i}}(\alpha_i)|\ n_{1},\dots,n_{d}\in\mathbb{N}\}$ where $k,r\in\mathbb{Z}_{+}$ and $\alpha_0,\alpha_1,\dots,\alpha_{d}\in G(K)$. Denote $\Phi=F^{k}$. Then $\Phi$ admits an equation $P(\Phi)=0$ in which $P(x)\in\mathbb{Z}[x]$ is a monic polynomial satisfying that every integer root of $P$ has the form $\pm p^{e}$ for some positive integer $e$. Let us analyze what the condition $S=\{\alpha_{0}+\sum\limits_{i=1}^{d} \Phi^{n_{i}}(\alpha_{i})|\ n_{1},\dots,n_{d}\in\mathbb{N}\}\subseteq\Gamma$ means. One can see that $\alpha_{0}+\alpha_{1}+\cdots+\alpha_{d}\in\Gamma$ and $\Phi(\alpha_{i})-\alpha_{i}\in\Gamma$ for each $1\leq i\leq d$. Write $\alpha_{0}+\alpha_{1}+\cdots+\alpha_{d}=c\cdot g$ and $\Phi(\alpha_{i})-\alpha_{i}=l_{i}\cdot g$ for some $c,l_{1},\dots,l_{d}\in\mathbb{Z}$. We may assume $l_{i}\neq0$ for each $i$ since otherwise $\Phi^{n_{i}}(\alpha_{i})=\alpha_{i}$ for all $n_{i}\in\mathbb{N}$ and so that this term can be absorbed into the constant term $\alpha_{0}$. Now the condition $S\subseteq\Gamma$ can be read as $l_{i}\Phi^{n}(g)\in\Gamma$ for any $1\leq i\leq d$ and $n\in\mathbb{N}$.

Write $l_{i}\Phi^{n}(g)=l_{i}^{(n)}g$ for some $l_{i}^{(n)}\in\mathbb{Z}$ (so $l_{i}^{(0)}=l_{i}$). Then $l_{i}^{(n)}g=l_{i}\Phi^{n}(g)=\Phi(l_{i}\Phi^{n-1}(g))=\Phi(l_{i}^{(n-1)}g)=l_{i}^{(n-1)}\Phi(g)$ for any $1\leq i\leq d$ and $n\in\mathbb{Z}_{+}$. Since $g$ is not a torsion point, we deduce $(l_{i}^{(n-1)})^{2}=l_{i}^{(n-2)}l_{i}^{(n)}$ for any $n\geq2$. Thus $l_{i}|l_{i}^{(1)}$ must hold because $l_{i}^{(0)}=l_{i}\neq0$ and each $l_{i}^{(n)}$ is an integer. So we can write $l_{i}^{(1)}=l_{i}t_{i}$ for some integer $t_{i}$ and hence $l_{i}(\Phi(g)-t_{i}g)=0$ for each $i$.

Since each $\Phi(g)-t_{i}g$ is a torsion point but $g$ is non-torsion, we can see that $t_{1}=\cdots=t_{d}:=t\in\mathbb{Z}$. Write $\Phi(g)=tg+h$ for some torsion point $h\in G(K)$. Then $\Phi^{n}(g)=t^{n}g+t^{n-1}h+t^{n-2}F(h)+\cdots+F^{n-1}(h)$ so that $\Phi^{n}(g)-t^{n}g$ is a torsion point for each $n\in\mathbb{Z}_{+}$. But $P(\Phi)=0$, so $P(t)\cdot g$ is a torsion point and hence $P(t)=0$ because $g$ is non-torsion. As a result, $t$ has the form $\pm p^{e}$ for some positive integer $e$.

Now we calculate $\Phi^{n_{i}}(\alpha_{i})$. We have $\Phi^{n_{i}}(\alpha_{i})=\alpha_{i}+(\Phi(\alpha_{i})-\alpha_{i})+\cdots+(\Phi^{n_{i}}(\alpha_{i})-\Phi^{n_{i}-1}(\alpha_{i}))=\alpha_{i}+l_{i}g+\cdots+l_{i}\Phi^{n_{i}-1}(g)=\alpha_{i}+l_{i}^{(0)}g+\cdots+l_{i}^{(n_{i}-1)}g$. Since $l_{i}^{(0)}=l_{i},l_{i}^{(1)}=l_{i}t$ and $(l_{i}^{(n-1)})^{2}=l_{i}^{(n-2)}l_{i}^{(n)}$ for any $n\geq2$, we deduce $l_{i}^{(n)}=l_{i}t^{n}$ for each $n\in\mathbb{N}$ by $l_{i}t\neq0$. Thus
$$\Phi^{n_{i}}(\alpha_{i})=\alpha_{i}+l_{i}(1+t+\cdots+t^{n_{i}-1})g=\alpha_{i}+l_{i}\frac{t^{n_{i}}-1}{t-1}\cdot g.$$

So we know
$$\alpha_{0}+\sum\limits_{i=1}^{d} \Phi^{n_{i}}(\alpha_{i})=(\alpha_{0}+\alpha_{1}+\cdots+\alpha_{d})+\sum\limits_{i=1}^{d} l_{i}\frac{t^{n_{i}}-1}{t-1}\cdot g=(c+\sum\limits_{i=1}^{d} l_{i}\frac{t^{n_{i}}-1}{t-1})\cdot g$$
for nonnegative integers $n_{1},\dots,n_{d}$ in which $c,l_{1},\dots,l_{d}\in\mathbb{Z}$ and $t=\pm p^{e}$ for some positive integer $e$. If $t$ is positive, we immediately see that\\
$$\{n\in\mathbb{Z}|\ n\cdot g\in S\}=\{c+\sum\limits_{i=1}^{d} l_{i}\frac{t^{n_{i}}-1}{t-1}|\ n_1,\dots,n_d\in\mathbb{N}\}=S_{p^{e},d,0}(c-\sum\limits_{i=1}^{d}\frac{l_i}{p^{e}-1};\frac{l_1}{p^{e}-1},\dots,\frac{l_d}{p^{e}-1}).$$
If $t$ is negative, we write $\{n\in\mathbb{Z}|\ n\cdot g\in S\}=\bigcup\limits_{\epsilon_{1},\cdots,\epsilon_{d}\in\{0,1\}} \{c+\sum\limits_{i=1}^{d} l_{i}\frac{(t+1)(t^{\epsilon_{i}}(t^{2})^{n_{i}}-1)}{t^{2}-1}|\ n_{1},\dots,n_{d}\in\mathbb{N}\}$ and then also conclude that it is a finite union of sets of the form $S_{p^{2e},d,0}(c_{0};c_{i0})$. Thus we are done.
\end{prf}

Now we have to consider the case that the infinite cyclic group $\Gamma\subseteq G(K)$ has no nontrivial $F^{k}$-invariant subgroup for any $k\in\mathbb{Z}_{+}$. Our result is as follows.

\begin{proposition}
Let $G$ be a semiabelian variety over an algebraically closed field $K$ of characteristic $p>0$ which is defined over a finite field $\mathbb{F}_{q}$. Let $F=\mathrm{Frob}_{q}$ be the Frobenius endomorphism of $G$ and let $\Gamma$ be the cyclic group generated by a non-torsion point $g\in G(K)$. Let $S$ be a subset of $\Gamma$ of the form $\{\alpha_0+\sum\limits_{i=1}^{d}\sum\limits_{j=0}^{r} F^{2^{j}n_i}(\alpha_{ij})|\ n_1,\dots,n_d\in\mathbb{N}\}$ where $d\in\mathbb{Z}_{+},r\in\mathbb{N}$ and $\alpha_{0},\alpha_{ij}\in G(K)$ where $(i,j)\in\{1,\dots,d\}\times\{0,\dots,r\}$. Suppose that
\begin{enumerate}
\item
$F$ admits an equation $P(F)=0$ where $P(x)=(x-q)P_{0}(x)$ and $P_0(x)\in\mathbb{Z}[x]$ is a monic polynomial satisfying that every complex root of $P_0$ has absolute value $q^{\frac{1}{2}}$, and

\item
$\Gamma$ has no nontrivial $F^{k}$-invariant subgroup for any $k\in\mathbb{Z}_{+}$.
\end{enumerate}
Then $\{n\in\mathbb{Z}|\ n\cdot g\in S\}$ is a set of the form $S_{q,d,r-1}(c_0;c_{ij})$ defined in Definition 1.1 where $c_0,c_{ij}\in\mathbb{Q}$ for $(i,j)\in\{1,\dots,d\}\times\{0,\dots,r-1\}$.
\end{proposition}

\begin{remark}
\begin{enumerate}
\item
The condition (i) for the Frobenius endomorphism $F$ is easy to fulfill. One just raise to a bigger finite field $\mathbb{F}_{q^k}$ on which the torus part of $G$ splits and then this condition holds for $\mathrm{Frob}_{q^{k}}$.

\item
Notice that in the conclusion, the power of $2$ is only up to $2^{r-1}$ while we have $2^{r}$ in the widely $F$-set.
\end{enumerate}
\end{remark}

In the following discussions, we will always inherit the setting as in Proposition 2.12. More precisely, we will always let $G,F,\Gamma,g$ be objects as in Proposition 2.12 satisfying those two conditions. We start with several lemmas.

\begin{lemma}
Let $\alpha\in G(K)$ be a point, $(l_n)_{n\in\mathbb{N}}$ be a sequence of integers and $(t_n)_{n\in\mathbb{N}}$ be a sequence of torsion points in $G(K)$. If $F^{n}(\alpha)=l_{n}\cdot g+t_n$ for every nonnegative integer $n$, then $l_n=0$ for every $n\in\mathbb{N}$ and hence $\alpha$ is a torsion point.
\end{lemma}

\begin{prf}
For every $n\geq0$, we have $l_{n+1}\cdot g+t_{n+1}=F^{n+1}(\alpha)=l_{n}F(g)+F(t_{n})$. So $l_{n+1}\cdot g=l_{n}F(g)+(F(t_{n})-t_{n+1})$ and hence $(l_{n+1}^{2}-l_{n}l_{n+2})\cdot g$ is a torsion point for every $n\geq0$. But $g$ is non-torsion and hence we conclude that $l_{n+1}^{2}=l_{n}l_{n+2}$ for each $n$.

If $l_{0}\neq0$, then we may deduce a contradiction as follows. Since $(l_n)_{n\in\mathbb{N}}$ is a sequence of integers, we deduce that there exists an integer $a$ such that $l_{n}=l_{0}a^{n}$ for every $n\geq0$ (one may consider the cases that whether $l_1=0$ and find that the conclusion always holds). So we have $\alpha=l_{0}\cdot g+t_0$ and $F(\alpha)=l_{1}\cdot g+t_1=l_{0}a\cdot g+t_1$ and therefore, $F(l_{0}g)=F(\alpha)-F(t_0)=al_{0}g+(t_1-F(t_0))$. But there exists a positive integer $D$ such that $D(t_1-F(t_0))=0$, so we have $F(Dl_0g)=a\cdot Dl_0g$ and thus the nontrivial subgroup of $\Gamma$ generated by $Dl_0\cdot g$ is $F$-invariant. So we get a contradiction and conclude that $l_0=0$.

Now as $l_0=0$, we have $\alpha$ is a torsion point. Therefore, $F^{n}(\alpha)$ and hence $l_{n}\cdot g$ is torsion for each $n$. So $l_n=0$ for each $n$ and thus we are done.
\end{prf}

\begin{lemma}
Let $r$ be a nonnegative integer, $\alpha_0,\dots,\alpha_r\in G(K)$ be closed points and let $(l_n)_{n\in\mathbb{N}}$ be a sequence of integers, $(t_n)_{n\in\mathbb{N}}$ be a sequence of torsion points in $G(K)$. If $\sum\limits_{j=0}^{r}F^{2^{j}n}(\alpha_{j})=l_{n}\cdot g+t_{n}$ for every $n\geq0$, then there exists $c_0,\dots,c_{r-1}\in\mathbb{Q}$ such that $l_{n}=\sum\limits_{j=0}^{r-1}c_{j}q^{2^{j}n}$ for every $n\geq0$ (where $q$ is the power of $p$ which corresponds to $F$, i.e. $F=\mathrm{Frob}_{q}$).
\end{lemma}

\begin{prf}
We prove by induction on $r$. We notice that although the data $G,\Gamma,g$ are fixed, we may change $q,F$ into a suitable power in the induction procedure. The case when $r=0$ is exactly Lemma 2.14 as in this case we have proved that $l_n=0$ for each $n$. Suppose the result holds for $r$, now we consider the case of $r+1$.

Recall we have assumed that $F$ admits an equation $P(F)=0$ where $P(x)=(x-q)P_{0}(x)$ and $P_0(x)\in\mathbb{Z}[x]$ is a monic polynomial satisfying that every complex root of $P_0$ has absolute value $q^{\frac{1}{2}}$. Write $P_{0}(x)=(x-z_1)\cdots(x-z_{s})$ for some complex numbers $z_1,\dots,z_s$ of absolute value $q^{\frac{1}{2}}$. For every $m\geq0$, we let $Q_{m}(x)=(x-q^{m})(x-z_{1}^{m})\cdots(x-z_{s}^{m})$ which is a $\mathbb{Z}$-coefficient polynomial and let $R(x)=Q_{2}(x)Q_{4}(x)\cdots Q_{2^{r+1}}(x)$. Then we have $Q_{m}(F^{m})=0$ for every $m$ and thus $R(F^{2})=R(F^{4})=\cdots=R(F^{2^{r+1}})=0$. So by doing some $\mathbb{Z}$-coefficient linear combinations, we can see that there are a sequence of integers $(l_{n}')_{n\in\mathbb{N}}$ and a sequence of torsion points $(t_{n}')_{n\in\mathbb{N}}$ such that $F^{n}(R(F)(\alpha_0))=l_{n}'\cdot g+t_{n}'$ for every $n\geq0$. Hence by Lemma 2.14, we conclude that $R(F)(\alpha_0)$ is a torsion point.

Notice that the greatest common divisor of $P(x)$ and $R(x)$ in $\mathbb{Q}[x]$ is 1 or $x-q$ since the only common root of them is $q$ (if there exists one) and the root $q$ has multiplicity 1 in $P(x)$. But now $R(F)(\alpha_0)$ is a torsion point and $P(F)(\alpha_0)=0$, so we can see that $(F-q)(\alpha_0)$ must be a torsion point. Now denote $l_{n}'=l_{n+1}-ql_{n}$ and $\alpha_{j}'=(F^{2^{j}}-q)(\alpha_{j})$ for every $n\geq0$ and every $1\leq j\leq r+1$, then we can see that there is a certian sequence of torsion points $(t_{n}')_{n\in\mathbb{N}}$ such that $\sum\limits_{j=1}^{r+1}F^{2^{j}n}(\alpha_{j}')=l_{n}'\cdot g+t_{n}'$ for every $n$. We want to use the induction hypothesis towards $q^{2},F^{2}=\mathrm{Frob}_{q^2}$ and $\Gamma=\mathbb{Z}\cdot g$.

We need to check that the two conditions in Proposition 2.12 are still valid: firstly, condition (i) still holds because $F^{2}=\mathrm{Frob}_{q^2}$ satisfies $Q_{2}(F^{2})=0$ where $Q_{2}(x)=(x-q^{2})(x-z_{1}^{2})\cdots(x-z_{s}^{2})$ and $(x-z_{1}^{2})\cdots(x-z_{s}^{2})$ is a monic $\mathbb{Z}$-coefficient polynomial whose complex roots have absolute value $q$; and secondly, condition (ii) tautologically holds as $\Gamma$ has no nontrivial $F^{2k}$-invariant subgroups for any $k\in\mathbb{Z}_+$. So now we can use the induction hypothesis to conclude that there exists $c_{1}'\dots,c_{r}'\in\mathbb{Q}$ such that $l_{n}'=\sum\limits_{j=1}^{r}c_{j}'q^{2^{j}n}$ for every $n\geq0$. As a result, for every $n\geq0$, we can calculate that
\begin{equation*}
\begin{split}
l_{n}=q^{n}l_{0}+\sum\limits_{m=0}^{n-1}q^{n-1-m}l_{m}'=q^{n}l_{0}+q^{n-1}\sum\limits_{j=1}^{r}(c_{j}'\sum\limits_{m=0}^{n-1}q^{(2^{j}-1)m})=q^{n}l_{0}+q^{n-1}\sum\limits_{j=1}^{r}c_{j}'\frac{q^{(2^{j}-1)n}-1}{q^{2^{j}-1}-1}\\=l_{0}q^{n}+\sum\limits_{j=1}^{r}c_{j}'\frac{q^{2^{j}n}-q^n}{q^{2^{j}}-q}=(l_0-\sum\limits_{j=1}^{r}\frac{c_{j}'}{q^{2^{j}}-q})q^{n}+\sum\limits_{j=1}^{r}\frac{c_{j}'}{q^{2^{j}}-q}q^{2^{j}n}=\sum\limits_{j=0}^{r}c_{j}q^{2^{j}n}
\end{split}
\end{equation*}
for some $c_0,\dots,c_r\in\mathbb{Q}$. Thus we have proved that the result also holds in the case of $r+1$. So we finish the proof by induction.
\end{prf}

Now we can prove Proposition 2.12.

\proof[Proof of Proposition 2.12]
Firstly, the condition $S\subseteq\Gamma$ implies that for every $1\leq i\leq d$, there exists a sequence of integers $(l_{in})_{n\in\mathbb{N}}$ such that $\sum\limits_{j=0}^{r}F^{2^{j}n}\circ(F^{2^{j}}-1)(\alpha_{ij})=l_{in}\cdot g$ for every $n\geq0$. So by Lemma 2.15, we know that there exists $c_{i,0}\dots,c_{i,r-1}\in\mathbb{Q}$ such that $l_{in}=\sum\limits_{j=0}^{r-1}c_{i,j}q^{2^{j}n}$ for every $(i,n)\in\{1,\dots,d\}\times\mathbb{N}$. Let $\alpha_{0}+\sum\limits_{i=1}^{d}\sum\limits_{j=0}^{r}\alpha_{ij}=l_0\cdot g$ for some $l_0\in\mathbb{Z}$. Then we can calculate that\\
\begin{equation*}
\begin{split}
\alpha_{0}+\sum\limits_{i=1}^{d}\sum\limits_{j=0}^{r}F^{2^{j}n_{i}}(\alpha_{ij})=(\alpha_{0}+\sum\limits_{i=1}^{d}\sum\limits_{j=0}^{r}\alpha_{ij})+\sum\limits_{i=1}^{d}\sum\limits_{j=0}^{r}\sum\limits_{m=0}^{n_{i}-1}F^{2^{j}m}\circ(F^{2^{j}}-1)(\alpha_{ij})\\=l_0\cdot g+\sum\limits_{i=1}^{d}\sum\limits_{m=0}^{n_{i}-1}l_{im}\cdot g=(l_0+\sum\limits_{i=1}^{d}\sum\limits_{m=0}^{n_{i}-1}\sum\limits_{j=0}^{r-1}c_{i,j}q^{2^{j}m})\cdot g=(l_0+\sum\limits_{i=1}^{d}\sum\limits_{j=0}^{r-1}c_{i,j}\frac{q^{2^{j}n_{i}}-1}{q^{2^{j}}-1})\cdot g
\end{split}
\end{equation*}
for every $(n_1,\dots,n_d)\in\mathbb{N}^{d}$. As a result, we have\\
$$\{n\in\mathbb{Z}|\ n\cdot g\in S\}=\{l_0+\sum\limits_{i=1}^{d}\sum\limits_{j=0}^{r-1}c_{i,j}\frac{q^{2^{j}n_{i}}-1}{q^{2^{j}}-1}|\ n_1,\dots,n_d\in\mathbb{N}\}=S_{q,d,r-1}(l_0-\sum\limits_{i=1}^{d}\sum\limits_{j=0}^{r-1}\frac{c_{i,j}}{q^{2^{j}}-1};\frac{c_{i,j}}{q^{2^{j}}-1}).$$

Thus we are done.
\endproof

\section{Translation of algebraic groups}
In this section, we prove that the statement of Theorem 1.5 holds for translation of algebraic groups. This special case of Theorem 1.5 is interesting as it connects the $p$-Mordell--Lang problem and the dynamical $p$-Mordell--Lang problem. The main result of this section is as follows.

\begin{theorem}
Let $G$ be an algebraic group over an algebraically closed field $K$ of characteristic $p>0$. Let $g\in G(K)$ be a closed point and let $X\subseteq G$ be a closed subvariety. Then $\{n\in\mathbb{Z}|\ g^{n}\in X(K)\}$ is a widely $p$-normal set as in Definition 1.1.
\end{theorem}

The proof of Theorem 3.1 has 4 steps: firstly, we prove that it holds for isotrivial semiabelian varieties; secondly, we prove that it holds for semiabelian varieties; thirdly, we prove that it holds for commutative group varieties; and lastly, we prove that it holds for arbitrary algebraic groups. We remark that in the first three steps, we will use the additive notation since the ambient algebraic group is commutative. For instance, we will write ``$n\cdot g$" instead of ``$g^{n}$" for a multiple of an element in $G(K)$. But we will go back to the multiplicative notation in the last step.

Now we start with the first step.

\begin{proposition}
Theorem 3.1 holds when $G$ is an isotrivial semiabelian variety, i.e. when $G$ is a semiabelian variety defined over $\overline{\mathbb{F}_{p}}$.
\end{proposition}

\begin{prf}
We prove by induction on $\mathrm{dim}(X)$. Notice that we only need to deal with the case when $X$ is an irreducible closed subvariety of $G$. The case when $\mathrm{dim}(X)=0$ is easy.

Firstly, we argue that we may assume $\mathrm{Stab}_{G}(X)=\{0\}$ without loss of generality. Here we denote $\mathrm{Stab}_{G}(X)$ as the smooth algebraic subgroup of $G$ satisfies $\mathrm{Stab}_{G}(X)(K)=\{a\in G(K)|\ a+X=X\}$. We know that $\mathrm{Stab}_{G}(X)$ (and in fact every smooth algebraic subgroup of $G$) is defined over $\overline{\mathbb{F}_{p}}$ as an algebraic subgroup of $G$. So $G/\mathrm{Stab}_{G}(X)$ is also an isotrivial semiabelian variety and now $X/\mathrm{Stab}_{G}(X)\subseteq G/\mathrm{Stab}_{G}(X)$ is an irreducible closed subvariety with trivial translation stabilizer. But if we denote $\pi:G\rightarrow G/\mathrm{Stab}_{G}(X)$ as the quotient map, we find that $\{n\in\mathbb{Z}|\ n\cdot g\in X(K)\}=\{n\in\mathbb{Z}|\ n\cdot\pi(g)\in(X/\mathrm{Stab}_{G}(X))(K)\}$. So we may assume that $\mathrm{Stab}_{G}(X)=\{0\}$ without loss of generality.

Now we find $q$ which is a power of $p$ such that $G$ is defined over $\mathbb{F}_{q}$ and the corresponding Frobenius endomorphism $F=\mathrm{Frob}_{q}$ satisfies the condition (i) in Proposition 2.12 (we can do this, see Remark 2.13(i)). We may assume $g\in G(K)$ is non-torsion since otherwise $\{n\in\mathbb{Z}|\ n\cdot g\in X(K)\}$ is a finite union of arithmetic progressions. We will deal with two cases that whether $\Gamma=\mathbb{Z}\cdot g$ has a nontrivial $F^{k}$-invariant subgroup for some positive integer $k$.

~

\textbf{Case 1:} There exist $k_0,n_0\in\mathbb{Z}_{+}$ such that $n_0\Gamma$ is an $F^{k_0}$-invariant group.

Now since $\{n\in\mathbb{Z}|\ n\cdot g\in X(K)\}=\bigcup\limits_{a=0}^{n_0-1}(a+n_0\cdot\{n\in\mathbb{Z}|\ (a+n_0n)\cdot g\in X(K)\})$, we only need to show that $\{n\in\mathbb{Z}|\ (a+n_0n)\cdot g\in X(K)\}$ is a widely $p$-normal set for every $0\leq a\leq n_0-1$ by Remark 1.2(ii). Write $\{n\in\mathbb{Z}|\ (a+n_0n)\cdot g\in X(K)\}=\{n\in\mathbb{Z}|\ n\cdot n_0g\in(-ag+X)(K)\}$ and then we can apply Theorem 2.4 for $F^{k_0}=\mathrm{Frob}_{q^{k_0}}$, the $F^{k_0}$-invariant group $n_{0}\Gamma$ and the closed subvariety $-ag+X\subseteq G$. Then we can conclude that $(-ag+X)(K)\cap n_0\Gamma$ is an $F^{k_0}$-normal set in $n_0\Gamma$. So $\{n\in\mathbb{Z}|\ n\cdot n_0g\in(-ag+X)(K)\}$ is a widely $p$-normal set by Proposition 2.11 for every $0\leq a\leq n_0-1$ and in fact, the ``widely" $p$-sets will not appear in here. So we have proved that $\{n\in\mathbb{Z}|\ n\cdot g\in X(K)\}$ is a widely $p$-normal set in this case.

~

\textbf{Case 2:} $\Gamma$ has no nontrivial $F^{k}$-invariant subgroup for any positive integer $k$.

Now let $\Gamma_0=\mathbb{Z}[F]\cdot g$ be a finitely generated $F$-invariant subgroup of $G(K)$. Using Theorem 2.4, we can write $X(K)\cap\Gamma_0=\bigcup\limits_{i=1}^{A}(S_i+\Lambda_i)\cup\bigcup\limits_{i=1}^{B}S_{i}'$ where $S_1,\dots,S_A,S_1',\dots,S_B'$ are $F$-sets in $\Gamma_0$ and $\Lambda_1,\dots,\Lambda_A\subseteq\Gamma_0$ are infinite subgroups since one can split the term whose $\Lambda$ is a finite group into a finite union of $F$-sets in $\Gamma_0$. Notice that $X(K)\cap\Gamma=\bigcup\limits_{i=1}^{A}(\overline{S_i+\Lambda_i}(K)\cap\Gamma)\cup\bigcup\limits_{i=1}^{B}(S_{i}'\cap\Gamma)$, we have $\{n\in\mathbb{Z}|\ n\cdot g\in X(K)\}=\bigcup\limits_{i=1}^{A}\{n\in\mathbb{Z}|\ n\cdot g\in\overline{S_i+\Lambda_i}(K)\}\cup\bigcup\limits_{i=1}^{B}\{n\in\mathbb{Z}|\ n\cdot g\in S_{i}'\cap\Gamma\}$.

For every $1\leq i\leq A$, we have $\overline{S_i+\Lambda_i}\subsetneqq X$ because it is contained in $X$ and it has a nontrivial translation stabilizer. So $\mathrm{dim}(\overline{S_i+\Lambda_i})<\mathrm{dim}(X)$ because we have assumed that $X$ is irreducible and hence we conclude that $\{n\in\mathbb{Z}|\ n\cdot g\in\overline{S_i+\Lambda_i}(K)\}$ is a widely $p$-normal set by induction hypothesis.

For every $1\leq i\leq B$, we have $S_{i}'\cap\Gamma$ is a finite union of widely $F$-sets by Proposition 2.6. Now since we have known that the two conditions in Proposition 2.12 holds for $F=\mathrm{Frob}_{q}$, we can see that they also holds for every power $F^{k}=\mathrm{Frob}_{q^k}$. So using Proposition 2.12, we conclude that $\{n\in\mathbb{Z}|\ n\cdot g\in S_{i}'\cap\Gamma\}$ is a widely $p$-normal set for every $1\leq i\leq B$.

Combining the arguments above, we deduce that in this case $\{n\in\mathbb{Z}|\ n\cdot g\in X(K)\}$ is a widely $p$-normal set.

~

All in all, we proved that $\{n\in\mathbb{Z}|\ n\cdot g\in X(K)\}$ is a widely $p$-normal set by induction.
\end{prf}

Next, we prove that Theorem 3.1 holds for arbitrary semiabelian varieties.

\begin{lemma}
Theorem 3.1 holds when $G$ is a semiabelian variety.
\end{lemma}

\begin{prf}
Denote $\Gamma=\mathbb{Z}\cdot g$. Firstly, we may assume that $X$ is an irreducible closed subvariety of $G$ and $\overline{X(K)\cap\Gamma}=X$ without loss of generality. Then using Theorem 2.1, we can see that there exists

$\bullet$ a semiabelian subvariety $G_1\subseteq G$,

$\bullet$ an isotrivial semiabelian variety $G_0$ over $K$,

$\bullet$ a closed subvariety $X_0\subseteq G_0$,

$\bullet$ a surjective algebraic group homomorphism $f:G_1\rightarrow G_0$, and

$\bullet$ a point $x_0\in G(K)$

such that $X=x_0+f^{-1}(X_0)$ as a closed subset of $G$.

Now we have $\{n\in\mathbb{Z}|\ n\cdot g\in X(K)\}=\{n\in\mathbb{Z}|\ ng-x_0\in G_1(K)\text{ and }f(ng-x_0)\in X_0(K)\}$. We may assume that the set $\{n\in\mathbb{Z}|\ ng-x_0\in G_1(K)\}$ is infinite since otherwise $\{n\in\mathbb{Z}|\ n\cdot g\in X(K)\}$ will be a finite set and then we are done. Then $\{n\in\mathbb{Z}|\ ng-x_0\in G_1(K)\}$ must be an infinite arithmetic progressions. So we may write $\{n\in\mathbb{Z}|\ ng-x_0\in G_1(K)\}=a+n_0\mathbb{Z}$ where $a$ is an integer which satisfies $ag-x_0\in G_1(K)$ and $n_0$ is a positive integer which satisfies $n_0g\in G_1(K)$. Then $\{n\in\mathbb{Z}|\ n\cdot g\in X(K)\}=a+n_0\cdot\{n\in\mathbb{Z}|\ f(ag-x_0)+n\cdot f(n_0g)\in X_0(K)\}=a+n_0\cdot\{n\in\mathbb{Z}|\ n\cdot f(n_0g)\in(-f(ag-x_0)+X_0)(K)\}$. So we reduce to the case of isotrivial semiabelian varieties and finish the proof by Proposition 3.2 and Remark 1.2(ii).
\end{prf}

\begin{remark}
At this point, we remark that the original version of $p$DML conjecture holds for translation of abelian varieties. More precisely, the set $\{n\in\mathbb{Z}|\ n\cdot g\in X(K)\}$ will be a $p$-normal set as in Definition 5.1 if the ambient algebraic group is an abelian variety. This statement can be proved by running the procedure of the proof once again for abelian varieties. The key point is that every complex root of the minimal polynomial of the Frobenius endomorphism of an isotrivial abelian variety has the same absolute value. Therefore, in the setting of Proposition 2.6 we have that the intersection set of an $F$-set with a cyclic group is also a finite union of $F$-sets (instead of widely $F$-sets).
\end{remark}

Next, we prove that Theorem 3.1 holds for commutative group varieties. We remark that the notion of group varieties here is a synonym of smooth algebraic groups, i.e. they do not need to be connected.

\begin{lemma}
Theorem 3.1 holds when $G$ is a commutative group variety.
\end{lemma}

\begin{prf}
Firstly, by $\cite[\text{Theorem}\ 5.6.3\text{(i)}]{Bri17}$, we know that $G$ admits an exact sequence $0\rightarrow S\rightarrow G\rightarrow Q\rightarrow 0$ where $S$ is a semiabelian variety over $K$ and $Q$ is commutative group variety over $K$ which has finite exponent. So in particular $Q(K)$ is a torsion group. As a result, there exists a positive integer $n_0$ such that $n_0\cdot g\in S(K)$. Then we may write $\{n\in\mathbb{Z}|\ n\cdot g\in X(K)\}=\bigcup\limits_{a=0}^{n_0-1}(a+n_0\cdot\{n\in\mathbb{Z}|\ n\cdot n_0g\in(-a+X)(K)\})$. But for every $0\leq a\leq n_0-1$, we can see that $\{n\in\mathbb{Z}|\ n\cdot n_0g\in(-a+X)(K)\}=\{n\in\mathbb{Z}|\ n\cdot n_0g\in((-a+X)\cap S)(K)\}$ is a widely $p$-normal set by Lemma 3.3. So we conclude that $\{n\in\mathbb{Z}|\ n\cdot g\in X(K)\}$ is a widely $p$-normal set by Remark 1.2(ii).
\end{prf}

Finally, we prove Theorem 3.1. Notice that we go back to the multiplicative notation here.

\proof[Proof of Theorem 3.1]
Firstly, one can see that there is a smooth commutative algebraic subgroup $H\subseteq G$ whose underlying set is $\overline{g^{\mathbb{Z}}}$. Denote $X_{0}=X\cap H$ which is a closed subvariety of $H$. Then $\{n\in\mathbb{Z}|\ g^{n}\in X(K)\}=\{n\in\mathbb{Z}|\ g^{n}\in X_0(K)\}$ is a widely $p$-normal set by Lemma 3.5 and thus we have finished the proof.
\endproof

At the end of this Section, we would like to give an example. It is well-known that the set $\{n\in\mathbb{Z}|\ g^{n}\in X(K)\}$ in Theorem 3.1 can be a ``$p$-set" when the ambient algebraic group is an algebraic torus. We shall give an explicit example to show that when the ambient algebraic group is an abelian variety, this set can also be something beyond a finite union of arithmetic progressions. The result in this example will be used in Section 5 to give a rigorous disproof of the original version of $p$DML conjecture.

\begin{example}
\begin{enumerate}
\item
Let $p=5$ and let $K=\overline{\mathbb{F}_{p}(t)}$. Let $E$ be the elliptic curve $x_{1}^{2}x_{2}=x_{0}^{3}+x_{2}^{3}$ in $\mathbb{P}_{K}^{2}$ with zero point $O=[0,1,0]\in E(K)$. Let $A=E\times E$ be an abelian variety. We embed $A$ into $\mathbb{P}_{K}^{8}$ by Segre embedding, i.e. $[x_0,x_1,x_2]\times[y_0,y_1,y_2]\mapsto[x_0y_0,x_0y_1,x_0y_2,x_1y_0,x_1y_1,x_1y_2,x_2y_0,$
$x_2y_1,x_2y_2]$. Let $z_{ij}$ be the coordinate of $\mathbb{P}^{8}$ corresponding to $x_iy_j$ for any $0\leq i,j\leq2$. Let $X\subseteq A$ be the closed subvariety $\{z_{02}=z_{20}+z_{22}\}\cap A$. Let $g=(Q_1,Q_2)\in A(K)$ where $Q_1=(t+1,\sqrt{(t+1)^{3}+1}),Q_2=(t,\sqrt{t^3+1})$ are points lie in the affine chart of $E(K)$. Denote $S=\{n\in\mathbb{N}|\ n\cdot g\in X(K)\}$. Then we have
\begin{enumerate}
\item
$\{p^{n}|\ n\in\mathbb{N}\}\subseteq S$, and
\item
$S\subseteq\{0\}\cup\{p^{k}m|\ k\in\mathbb{N},m\in\mathbb{Z}_{+},m\equiv\pm1\ (\mathrm{mod}\ 2p)\}$.
\end{enumerate}

Thus $S$ cannot be a finite union of arithmetic progressions in $\mathbb{N}$.

\item
Let $A$ be an abelian surface over an algebraically closed field $K$ of characteristic $p>0$. Let $X\subseteq A$ be a closed subvariety and let $g\in A(K)$ be a closed point. Then $\{n\in\mathbb{Z}|\ n\cdot g\in X(K)\}$ is a union of finitely many arithmetic progressions (possibly singleton) along with finitely many sets of the form $S_{q,1,0}(\frac{c_0}{q-1};\frac{c_1}{q-1})$ as in Definition 1.1 where $q$ is a power of $p$ and $c_0,c_1$ are integers satisfying $q-1\mid c_0+c_1$.
\end{enumerate}
\end{example}

\begin{prf}
\begin{enumerate}
\item
Let $F$ be the Frobenius endomorphism $\mathrm{Frob}_{p}$ of $E$. Since $E$ is a supersingular elliptic curve, we have $F^{2}=[-p]\in\mathrm{End}(E)$. As a result, $p^{n}\cdot P=(-1)^{n}\cdot F^{2n}(P)$ for any $n\in\mathbb{N}$ and $P\in E(K)$. Thus we can see that (i) holds.

To prove (ii), we need the explicit formula of the multiplication-by-$m$ map of an elliptic curve described in $\cite[(\text{III, Ex. 3.7})]{Sil09}$. We apply this result to our $E$.

For any positive integer $m$, there exist coprime polynomials $f_{m}(x)=x^{m^{2}}+$ (lower order terms) and $g_{m}(x)=m^{2}x^{m^{2}-1}+$ (lower order terms) in $\mathbb{F}_{p}[x]$ such that

$$m\cdot P=\left\{
\begin{array}{cc}
\left(\frac{f_{m}(x)}{g_{m}(x)},y_m\right), & g_{m}(x)\neq0 \\
O, & g_{m}(x)=0
\end{array}
\right.$$
where $P=(x,y)$ lies in the affine chart of $E(K)$ and $y_m$ is a certain element of $K$. In particular, the points $Q_1$ and $Q_2$ are non-torsion.

Now $S=\{0\}\cup\{m\in\mathbb{Z}_{+}|\ \frac{f_{m}(t+1)}{g_{m}(t+1)}=\frac{f_{m}(t)}{g_{m}(t)}+1\}$. But since $f_{m}(x)$ and $g_{m}(x)$ are coprime, this condition on $m$ yields $g_{m}(t+1)=g_{m}(t)$. So $g_{m}(x)$ must be a polynomial of $x^{p}-x$ and as a result, the number of different roots of $g_{m}(x)$ in $K$ is a multiple of $p$.

Denote this number by $pd_{m}$ and write $m=p^{k}m'$ in which $p\nmid m'$. Then by the supersingularity of $E$, we deduce
$$m'^{2}=|E[m]|=\left\{
\begin{array}{cc}
1+2pd_{m}, & 2\nmid m \\
4+2(pd_{m}-3), & 2|m
\end{array}
\right.$$
But $p=5$ cannot be a factor of $m'^{2}+2$, so we have $2\nmid m$ and $m'\equiv\pm1\ (\mathrm{mod}\ p)$. Thus we are done.

\item
Assume $g$ is non-torsion without loss of generality. We only have to deal with the case that $X\subseteq A$ is an irreducible curve and $\{n\in\mathbb{Z}|\ n\cdot g\in X(K)\}$ is an infinite set. We may further assume that $n_0g+X\neq X$ for every positive integer $n_0$ because otherwise $\{n\in\mathbb{Z}|\ n\cdot g\in X(K)\}$ will be a finite union of arithmetic progressions and hence we are done. Now recall the statement in Remark 3.4 which implies that $\{n\in\mathbb{Z}|\ n\cdot g\in X(K)\}$ is a union of finitely many arithmetic progressions along with finitely many sets of the form $S_{q,d,0}(\frac{c_0}{q-1};\frac{c_i}{q-1})$ where $c_0,c_1,\dots,c_d$ are integers satisfying $q-1\mid c_0+c_1+\cdots+c_d$.

Notice that for every positive integer $n_0$, the set $\{n\in\mathbb{Z}|\ ng\in X(K)\text{ and }(n+n_0)g\in X(K)\}$ must be finite because of the assumption $n_0g+X\neq X$. But then this condition forces $\{n\in\mathbb{Z}|\ n\cdot g\in X(K)\}$ to be a union of finitely many singletons along with finitely many sets of the form $S_{q,1,0}(\frac{c_0}{q-1};\frac{c_1}{q-1})$ where $c_0,c_1$ are integers satisfying $q-1\mid c_0+c_1$. Thus we have finished the proof.
\end{enumerate}
\end{prf}

\begin{remark}
We find it hard to completely determine the return set. In part (i) above, we do not know how to prove that $S=\{0\}\cup\{p^{n}|\ n\in\mathbb{N}\}$ although we believe this is true. Maybe one can follow the procedure in $\cite{BGM}$ to get a rigorous proof.
\end{remark}

\section{Bounded-degree self-maps}
We will prove Theorem 1.5 in this section. We shall use some knowledge in $\cite{Bri19}$ as well as a regularization theorem (Theorem 4.6) to deduce Theorem 1.5 from Theorem 3.1. More precisely, we will deal with the case of bounded-degree automorphisms in subsection 4.1, and prove Theorem 1.5 in subsection 4.2 by reducing to the case of bounded-degree automorphisms.

Throughout this section, we fix an algebraically closed field $K$ of characteristic $p>0$ and let everything be over this field. We require $K$ to be of positive characteristic only because we use Theorem 3.1 in the proof of Proposition 4.4.

\subsection{Bounded-degree automorphisms}

In this subsection, let $X$ be a projective variety and let $f$ be a bounded-degree automorphism of $X$. We denote $\mathrm{N}^{1}(X)$ as the group of line bundles on $X$ up to numerical equivalence, which is a finite free $\mathbb{Z}$-module. For $L\in\mathrm{Pic}(X)$, we denote by $[L]_{\mathrm{num}}$ the class of $L$ in $\mathrm{N}^{1}(X)$. We will prove that Theorem 1.5 holds for the bounded-degree system $(X,f)$ at the end of this subsection. Now we start with the following proposition which says that bounded-degree automorphisms come from group actions.

\begin{proposition}
Let $X,f$ be as above. Then there exists a (not necessarily connected) group variety $G$, a group action $F:G\times X\rightarrow X$ and a point $g_{0}\in G(K)$ such that $f=F_{g_{0}}$ in which $F_{g_{0}}$ is the automorphism of $X$ induced by the group action.
\end{proposition}

Firstly, we shall show that the action of the bounded-degree automorphism $f$ on $\mathrm{N}^{1}(X)$ is unipotent. We need the following lemma on intersection theory.

\begin{lemma}
Let $X$ be as above and let $C\subseteq X$ be an irreducible closed subcurve. Then there exists an ample line bundle $L_{0}\in\mathrm{Pic}(X)$ such that for every ample line bundle $L\in\mathrm{Pic}(X)$, we have $L\cdot L_{0}^{\mathrm{dim}X-1}\geq L\cdot C$.
\end{lemma}

\begin{prf}
We prove the assertion by induction on $\mathrm{dim}X$. If $\mathrm{dim}X=1$, then there is nothing to prove. So we may assume that $\mathrm{dim}X\geq2$.

Denote $I_{C}\subseteq\mathcal{O}_{X}$ as the ideal sheaf of $C\subseteq X$. Pick an ample line bundle $L_{1}\in\mathrm{Pic}(X)$ such that $I_{C}\otimes L_{1}$ is globally generated. Pick a nonzero global section $s\in\Gamma(X,I_{C}\otimes L_{1})\subseteq\Gamma(X,L_{1})$ (notice that $I_{C}\otimes L_{1}$ is a subsheaf of $L_{1}$). Let $D=(s)_{0}$ be the divisor of zeros of $s$, which is an effective Cartier divisor on $X$ such that $L_{1}\cong\mathscr{L}(D)$. Let $Y\subseteq X$ be the closed subscheme associated with $D$. Then $C\subseteq Y$ and $Y$ is of pure codimension 1 in $X$.

Let $Y_{0}\subseteq Y$ be an irreducible component of $Y$ containing $C$ and equip $Y_{0}$ with the reduced induced closed subscheme structure. Then $Y_{0}$ is a projective variety of dimension $\mathrm{dim}X-1$. Let $i:Y_{0}\hookrightarrow X$ be the closed immersion. By induction hypothesis, there is an ample line bundle $L_{2}\in\mathrm{Pic}(Y_{0})$ such that $(L'\cdot L_{2}^{\mathrm{dim}Y_{0}-1})_{Y_{0}}\geq(L'\cdot C)_{Y_{0}}$ for every ample line bundle $L'\in\mathrm{Pic}(Y_{0})$. Now we choose an ample line bundle $L_{0}\in\mathrm{Pic}(X)$ such that both $L_{0}-L_{1}$ and $i^{*}L_{0}-L_{2}$ are ample. We claim that $L_{0}$ has the desired property.

Indeed, for every ample line bundle $L\in\mathrm{Pic}(X)$, we have $(L\cdot L_{0}^{\mathrm{dim}X-1})\geq(L\cdot L_{0}^{\mathrm{dim}X-2}\cdot L_{1})=(L\cdot L_{0}^{\mathrm{dim}X-2}\cdot Y)\geq(L\cdot L_{0}^{\mathrm{dim}X-2}\cdot Y_{0})=(i^{*}L\cdot(i^{*}L_{0})^{\mathrm{dim}X-2})_{Y_{0}}\geq(i^{*}L\cdot(L_{2})^{\mathrm{dim}X-2})_{Y_{0}}\geq(i^{*}L\cdot C)_{Y_{0}}=(L\cdot C)$. Thus we finish the proof by induction.
\end{prf}

\begin{lemma}
Let $X,f$ be as above. Then there exists a positive integer $n_{0}$ such that $(f^{n_{0}})^{*}:\mathrm{N}^{1}(X)\rightarrow\mathrm{N}^{1}(X)$ is the identity map.
\end{lemma}

\begin{prf}
Pick a $\mathbb{Z}$-basis $\{[L_{1}]_{\mathrm{num}},\dots,[L_{d}]_{\mathrm{num}}\}$ of $\mathrm{N}^{1}(X)$. Then there exists $\{C_{1},\dots,C_{d}\}$ which are $\mathbb{Q}$-coefficient 1-cycles in $X$, such that $L_{i}\cdot C_{j}=\delta_{ij}$ for all $1\leq i,j\leq d$ where $\delta_{ij}$ is the Kronecker symbol. Let $A\in GL_{d}(\mathbb{Z})$ be the matrix corresponds to $f^{*}:\mathrm{N}^{1}(X)\rightarrow\mathrm{N}^{1}(X)$ under this basis. We have to show that there is a poistive integer $n_{0}$ such that $A^{n_{0}}=I_{d}$ which is the identity matrix.

For each nonnegative integer $n$, we have $A^{n}=((f^{n})^{*}(L_{1}),\dots,(f^{n})^{*}(L_{d}))^{\top}\cdot(C_{1},\dots,C_{d})$. But by Remark 1.4 and Lemma 4.2, we can see that each sequence $\{(f^{n})^{*}(L_{i})\cdot C_{j}|\ n\in\mathbb{N}\}$ is bounded. So $\{A^{n}|\ n\in\mathbb{N}\}$ is a sequence in $GL_{d}(\mathbb{Z})$ in which each element is bounded. As a result, there are only finitely many different matrices in that sequence. Thus we are done because $A$ is invertible.
\end{prf}

Next, we need to recall some knowledge of the automorphism groups of projective varieties. For a reference, see $\cite[\mathrm{Section}\ 2]{Bri19}$.

Let $X$ be our projective variety. Let $\mathbf{Aut}_{X}$ be the contravariant functor from the category of (locally noetherian) $K$-schemes to the category of groups, which sends the $K$-scheme $S$ to the group $\mathrm{Aut}(X\times S/S)$ (the products will always be taken over $K$). This functor is represented by a locally algebraic group $\mathrm{Aut}_{X}$ over $K$. Let $\mathrm{Aut}(X)=\mathrm{Aut}_{X,\mathrm{red}}$ be the reduced closed (locally algebraic) subgroup of $\mathrm{Aut}_{X}$, then there are canonical bijections $\mathrm{Aut}(X/K)=\mathrm{Aut}_{X}(K)=\mathrm{Aut}(X)(K)$. Let $\mathrm{Aut}^{0}(X)$ be the identity component of $\mathrm{Aut}(X)$, which is a connected group variety. Then $\mathrm{Aut}(X)$ acts on $\mathrm{N}^{1}(X)$ and in fact $\mathrm{Aut}^{0}(X)$ acts trivially on it (see $\cite[\mathrm{Lemma}\ 2.8]{Bri19}$ and the discussion above it).

Let $L\in\mathrm{Pic}(X)$ be an ample line bundle and let Aut($X,[L]_{\mathrm{num}}$) be the stabilizer of $[L]_{\mathrm{num}}$ under the action of $\mathrm{Aut}(X)$. Then $\cite[\mathrm{Theorem}\ 2.10]{Bri19}$ says that Aut($X,[L]_{\mathrm{num}}$) is a closed algebraic subgroup of $\mathrm{Aut}(X)$. Notice that $\mathrm{Aut}(X,[L]_{\mathrm{num}})(K)\subseteq\mathrm{Aut}(X)(K)$ is canonically identified with $\{f\in\mathrm{Aut}(X/K)|\ f^{*}(L)\equiv L\}\subseteq\mathrm{Aut}(X/K)$ where ``$\equiv$" stands for numerical equivalence.

Now we can prove Proposition 4.1.

\proof[Proof of Proposition 4.1]
Let $g_{0}\in\mathrm{Aut}_{X}(K)=\mathrm{Aut}(X)(K)$ be the closed point which corresponds to the bounded-degree automorphism $f$. Combining Lemma 4.3 and the discussion above, we can see that there is a positive integer $n_{0}$ such that $g_{0}^{n_{0}}$ lies in a closed algebraic subgroup of $\mathrm{Aut}(X)$ (as it lies in any $\mathrm{Aut}(X,[L]_{\mathrm{num}})$ where $L\in\mathrm{Pic}(X)$ is an ample line bundle). As a result, if we let $G\subseteq\mathrm{Aut}_{X}$ be the closed smooth (locally algebraic) subgroup whose underlying space is $\overline{\{g_{0}^{n}|\ n\in\mathbb{Z}\}}$, then $G$ is in fact an (algebraic) group variety.

Now notice that there is a natural group action $\sigma:\mathrm{Aut}_{X}\times X\rightarrow X$ such that $f=\sigma_{g_{0}}$ (in which $\sigma_{g_{0}}$ is the automorphism of $X$ induced by the group action), we may just let $F:G\times X\rightarrow X$ be the group action induced by $\sigma$ and then one can verify that the conclusion of Proposition 4.1 holds for $G,F$ and $g_0\in G(K)$.
\endproof

Now we can prove Theorem 1.5 for bounded-degree automorphisms of projective varieties.

\begin{proposition}
Let $X$ and $f$ be as in the beginning of this subsection. Then the conclusion of Theorem 1.5 holds for the dynamical system $(X,f)$.
\end{proposition}

\begin{prf}
Let $x\in X(K)$ be a closed point and let $V\subseteq X$ be a closed subvariety. We will prove that $\{n\in\mathbb{Z}|\ f^{n}(x)\in V(K)\}$ is a widely $p$-normal set.

Let $G,F$ and $g_0\in G(K)$ be as in Proposition 4.1. Since $f=F_{g_0}$, we know that $f^{n}=F_{g_{0}^{n}}$ for every integer $n$ where $F_{g_{0}^{n}}$ is the automorphism of $X$ induced by $F$ and $g_{0}^{n}\in G(K)$. So for each integer $n$, we have $f^{n}(x)=F(g_{0}^{n},x)$. Now let $i_{x}:G\hookrightarrow G\times X$ be the closed immersion given by $g\mapsto(g,x)$ and let $j:G\rightarrow X$ be the composition $F\circ i_{x}$. Then we have $\{n\in\mathbb{Z}|\ f^{n}(x)\in V(K)\}=\{n\in\mathbb{Z}|\ g_{0}^{n}\in j^{-1}(V)\}$. Thus the result follows from Theorem 3.1.
\end{prf}

\subsection{Proof of Theorem 1.5}

We will finish the proof of Theorem 1.5 in this subsection. Firstly, we reduce to the case in which the orbit $\mathcal{O}_{f}(x)$ is dense in $X$.

\begin{lemma}
In order to prove Theorem 1.5, we may assume that the orbit $\mathcal{O}_{f}(x)$ is dense in $X$ without loss of generality.
\end{lemma}

\begin{prf}
Assume that we have proved Theorem 1.5 with the additional assumption that $\mathcal{O}_{f}(x)$ is dense in $X$, we want to prove Theorem 1.5. We will do induction on $\mathrm{dim}(X)$. The case in which $\dim(X)=1$ is easy.

Now assume $f:X\dashrightarrow X$ is a bounded-degree self-map of a projective variety $X$ and $x\in X(K)$ is a point such that $\mathcal{O}_{f}(x)$ is well-defined. Let $V\subseteq X$ be a closed subvariety. We will prove that $\{n\in\mathbb{N}|\ f^{n}(x)\in V(K)\}$ is a widely $p$-normal set in $\mathbb{N}$. Assume further without loss of generality that $\mathcal{O}_{f}(x)$ is not dense in $X$. By substituting $x$ by a proper iterate, we may assume that the closed subsets $\overline{\{f^{n}(x)|\ n\geq N\}}$ are all the same for any nonnegative integer $N$. Denote this closed subset as a proper closed subvariety $X_{0}\subseteq X$.

Let $X_{11},\dots,X_{1d_{1}},X_{21},\dots,X_{2d_{2}},\dots,X_{r1},\dots,X_{rd_{r}}$ be the irreducible components of $X_{0}$ such that $\mathrm{dim}(X_{11})=\dots=\mathrm{dim}(X_{1d_{1}})>\mathrm{dim}(X_{21})=\dots=\mathrm{dim}(X_{2d_{2}})>\dots>\mathrm{dim}(X_{r1})=\dots=\mathrm{dim}(X_{rd_{r}})$. Notice that since $\mathcal{O}_{f}(x)$ is dense in $X_{0}$, we have $\mathcal{O}_{f}(x)\cap X_{ij}$ is nonempty for each $i,j$. So $U_{ij}:=\mathrm{Dom}(f)\cap X_{ij}$ is a nonempty open subset of $X_{ij}$ for each $i,j$. We can see that $\bigcup\limits_{i,j}f(U_{ij})$ is a dense subset of $X_{0}$.

As a result, one can choose a pair $(\sigma_{1}(i,j),\sigma_{2}(i,j))$ for each $i,j$ such that $f(U_{ij})\subseteq X_{\sigma_{1}(i,j)\sigma_{2}(i,j)}$. Using the density of $\bigcup\limits_{i,j}f(U_{ij})$ in $X_{0}$, one can show that (by induction on $i$ from 1 to $r$) $\sigma_{1}(i,j)=i$ for each $i,j$ and $\sigma_{2}(i,1),\dots,\sigma_{2}(i,d_{i})$ is a permutation of $1,\dots,d_{i}$ for each $i=1,2,\dots,r$. Furthermore, by the same reason, $f(U_{ij})$ must be dense in $X_{i\sigma_{2}(i,j)}$ for each $i,j$. We abbreviate $\sigma_{2}(i,j)$ as $\sigma(i,j)$.

By the discussion above, we see that $f$ induces dominant rational maps $f_{ij}:X_{ij}\dashrightarrow X_{i\sigma(i,j)}$ for each $i,j$, and $U_{ij}\subseteq\mathrm{Dom}(f_{ij})$. Suppose $x\in U_{i_{0}j_{1}}(K)$ and let $j_{1},j_{2},\dots,j_{t},j_{t+1}=j_{1}$ be a circle under the action of $\sigma(i_{0},j)$. We abbreviate $X_{i_{0}j_{k}}$ as $X_{k}$ and $f_{i_{0}j_{k}}$ as $f_{k}$, then we get a circle of dominant rational maps $f_{k}:X_{k}\dashrightarrow X_{k+1}$ ($k=1,\dots,t$, understood by modulo $t$ for the indices), all induced by $f$. Denote $g_{k}:X_{k}\dashrightarrow X_{k}$ as the composite $f_{k+t-1}\circ\cdots\circ f_{k}$ for $k=1,\dots,t$, which is a dominant rational self-map of $X_{k}$.

One can verify that for each $k=1,\dots,t$, $\mathrm{Dom}(f^{t})\cap X_{k}$ is nonempty and $f^{t}|_{X_{k}}$ maps into $X_{k}$. In fact, one can show that $f^{t}|_{X_{k}}:X_{k}\dashrightarrow X_{k}$ is same as $g_{k}$. So by (the proof of) $\cite[\mathrm{Proposition}\ 3.2]{JXSZ}$, we know that $g_{k}$ are bounded-degree (dominant) self-maps for each $k=1,\dots,t$. Moreover, since $x\in U_{i_{0}j_{1}}(K)\subseteq X_{1}(K)$, we have $f^{n}(x)=f_{n}\circ\cdots\circ f_{1}(x)\in U_{i_{0}j_{n+1}}(K)\subseteq X_{n+1}(K)$ for each $n\in\mathbb{N}$. As a result, we have $\mathcal{O}_{g_{k}}(f^{k-1}(x))$ is well-defined and $g_{k}^{n}(f^{k-1}(x))=f^{nt+k-1}(x)$ for every $k=1,\dots,t$ and every $n\in\mathbb{N}$. Notice that $\mathrm{dim}(X_{1})=\dots=\mathrm{dim}(X_{t})<\mathrm{dim}(X)$, we know each $g_{k}$ satisfies the conclusion of Theorem 1.5 in view of the induction hypothesis. So using Remark 1.2(ii), we may conclude that $\{n\in\mathbb{N}|\ f^{n}(x)\in V(K)\}=\bigcup\limits_{k=1}^{t}((k-1)+t\cdot\{n\in\mathbb{N}|\ f^{nt+k-1}(x)\in V(K)\})=\bigcup\limits_{k=1}^{t}((k-1)+t\cdot\{n\in\mathbb{N}|\ g_{k}^{n}(f^{k-1}(x))\in(V\cap X_k)(K)\})$ is a widely $p$-normal set in $\mathbb{N}$ for any closed subvariety $V\subseteq X$. Hence we finish the proof by induction.
\end{prf}

Next, we introduce our main tool. This result was firstly proved in $\cite[\text{Section}\ 5]{HZ96}$ for the case that the base field is $\mathbb{C}$. By adapting the methods in $\cite{HZ96}$, the second author wrote a note $\cite{Yang}$ which deals with the arbitrary characteristic case and copes with the details carefully. One can consult $\cite[\text{Section}\ 5]{HZ96},\cite[\text{Theorem}\ 2.5]{Can14}$ or $\cite[\mathrm{Corollary}\ 1.3]{Yang}$ for references.

\begin{theorem}
Let $f:X\dashrightarrow X$ be a bounded-degree self-map of a projective variety. Then there exists a projective variety $Y$, a birational map $\pi:Y\dashrightarrow X$ and a bounded-degree automorphism $g$ of $Y$, such that $f\circ\pi=\pi\circ g$.
\end{theorem}

We need another lemma before proving Theorem 1.5.

\begin{lemma}
Let $X$ and $Y$ be two varieties and let $\pi:Y\dashrightarrow X$ be a dominant rational map. Let $f$ be a dominant rational self-map of $X$ and $g$ be a flat endomorphism of $Y$ such that $f\circ\pi=\pi\circ g$. Let $x\in X(K)$ be a point such that the orbit $\mathcal{O}_{f}(x)$ is well-defined and dense in $X$. Suppose that $\{n\in\mathbb{N}|\ g^{n}(y)\in W(K)\}$ is a widely $p$-normal set in $\mathbb{N}$ for every point $g\in Y(K)$ and every closed subvariety $W\subseteq Y$. Then for every closed subvariety $V\subseteq X$, we have $\{n\in\mathbb{N}|\ f^{n}(x)\in V(K)\}$ is a widely $p$-normal set in $\mathbb{N}$.
\end{lemma}

\begin{prf}
Firstly, notice that since the flat endomorphism $g$ of $Y$ must be dominant, no problem will occur when compositing the maps. Let $D=\text{Dom}(\pi)$ and let $\pi_0:D\rightarrow X$ be the morphism which represents $\pi$. Then $D$ is an open dense subset of $Y$ and $\pi_0(D)$ is a constructible dense subset of $X$ since $\pi$ is dominant. So $\pi_0(D)$ contains an open dense subset of $X$. By substituting $x$ by a proper iterate, we may assume $x\in\pi_0(D)$ without loss of generality because $\mathcal{O}_{f}(x)$ is dense in $X$. Then we can choose a point $y\in D(K)$ such that $\pi_0(y)=x$. We firstly prove that $\mathcal{O}_{g}(y)\subseteq D$.

Let $n$ be a nonnegative integer. We want to prove that $g^{n}(y)\in D$. By the assumption that $\mathcal{O}_{f}(x)$ is well-defined, we can see that $x\in\text{Dom}(f^{n})$. As a result, we have $y\in\text{Dom}(f^{n}\circ\pi)$. So $y\in\text{Dom}(\pi\circ g^{n})$ because $f^{n}\circ\pi=\pi\circ g^{n}$. Now since $g^{n}$ is a flat morphism, we conclude that $g^{n}(y)\in D$ by $\cite[\text{Proposition}\ 20.3.11]{GD67}$. Therefore, we deduce that $g^{n}(y)\in D$ for each nonnegative integer $n$ and thus $\mathcal{O}_{g}(y)\subseteq D$.

Now for each nonnegative integer $n$, we have $f^{n}(x)=f^{n}(\pi_0(y))=\pi_0(g^{n}(y))$ since $x\in\text{Dom}(f^{n})$, $g^{n}(y)\in D$ and $f^{n}\circ\pi=\pi\circ g^{n}$. Let $W$ be the closure of $\pi_{0}^{-1}(V)$ in $Y$. Then $W$ is a closed subvariety of $Y$ such that $W\cap D=\pi_{0}^{-1}(V)$. Now combining the fact $\mathcal{O}_{g}(y)\subseteq D$ with the equality $f^{n}(x)=\pi_0(g^{n}(y))$ above, we know that $\{n\in\mathbb{N}|\ f^{n}(x)\in V(K)\}=\{n\in\mathbb{N}|\ g^{n}(y)\in W(K)\}$. Hence the result follows.
\end{prf}

\begin{remark}
The analogue of the lemma above for the 0-DML property (i.e. the statement which asserts that the return set is a finite union of arithmetic progressions in $\mathbb{N}$) is also valid. No change is needed to make in the proof.
\end{remark}

Now we can finish the proof of Theorem 1.5.

\proof[Proof of Theorem 1.5]
Combining Lemma 4.5, Theorem 4.6, Lemma 4.7 and Proposition 4.4 and then we are done.
\endproof

At the end of this section, we would like to give a remark on the complexity of the return set.

\begin{remark}
By running the proof procedure more carefully, one may let all of the ``widely $p$-sets" involved in the return set have the form $\{c_{0}+\sum\limits_{i=1}^{d}\sum\limits_{j=0}^{r_i} c_{ij}q^{2^{j}n_{i}}|\ n_{1},\dots,n_{d}\in\mathbb{N}\}$ where $q$ is a power of $p$ and $d+r_1+\cdots+r_d+|\{i|\ r_i>0,1\leq i\leq d\}|$ is not greater than the dimension of the closed subvariety $V$ at least in the situation of Theorem 3.1 or Proposition 4.4. The key point is that by running through the proof in $\cite{MS04}$, one may let all of the $F$-sets invloved in $X(K)\cap\Gamma$ have the form $\{\alpha_0+\sum\limits_{i=1}^{d}F^{kn_i}(\alpha_{i})|\ n_1,\dots,n_d\in\mathbb{N}\}$ where $d\leq\mathrm{dim}(X)$ in the situation of Theorem 2.4. Moreover, as we have seen in Proposotion 2.12 or Lemma 2.15, an extra place is needed when a real ``widely $p$-component" appears. Notice that the statement of Remark 2.10 is also needed here. This is because we have to make sure that for an $F$-set $\mathcal{F}$ and a cyclic group $C$, the ``length" of the widely $F$-sets appeared in the expression of $\mathcal{F}\cap C$ cannot be greater than the ``length" of $\mathcal{F}$.

We guess that with some more efforts, one may prove that this kind of statement also holds in the situation of Theorem 1.5.
\end{remark}

\section{A disproof of the original version of $p$DML conjecture}
In this section, we disprove the original version of the $p$DML conjecture. Firstly, let us review its statement.

\begin{definition}
Let $p$ be a prime. A $p$\emph{-normal set} (in $\mathbb{Z}$) is a union of finitely many arithmetic progressions (possibly singleton) along with finitely many sets of the form $S_{q,d,0}(\frac{c_{0}}{q-1};\frac{c_{i}}{q-1})$ (see Definition 1.1) in which $q$ is a power of $p$ and $c_0,c_1,\dots,c_d$ are integers satisfying $q-1\mid c_0+c_1+\cdots+c_d$. A $p$\emph{-normal set in} $\mathbb{N}$ is a subset of $\mathbb{N}$ which is, up to a finite set, equal to the intersection of a $p$-normal set and $\mathbb{N}$.
\end{definition}

As we have mentioned in the Introduction, this definition was firstly introduced in $\cite{Der07}$ for the Skolem--Mahler--Lech problem in positive characteristic. The original version of the $p$DML conjecture $\cite[\mathrm{Conjecture}\ 13.2.0.1]{BGT16}$ can be read as follows.

\begin{conjecture}
Let $X$ be a variety over an algebraically closed field $K$ of characteristic $p>0$ and let $f$ be a rational self-map of $X$. Let $x\in X(K)$ be a closed point such that the orbit $\mathcal{O}_{f}(x)$ is well-defined and let $V\subseteq X$ be a closed subvariety. Then $\{n\in\mathbb{N}|\ f^{n}(x)\in V(K)\}$ is a $p$-normal set in $\mathbb{N}$.
\end{conjecture}

We shall prove that this statement fails even in the special case of translations of isotrivial semiabelian varieties (but it is valid in the case of translation of abelian varieties, see Remark 3.4). We will give some heuristic of counterexamples in subsection 5.1 and then provide a rigorous disproof in subsection 5.2.

\subsection{Heuristic of counterexamples}
In this subsection, we provide some evidence to the existence of counterexamples of Conjecture 5.2 in the case of translation of isotrivial semiabelian varieties. Since we only want to give some heuristic arguments here, we will not be quite rigorous in this subsection.

\begin{example}
Let $G$ be a semiabelian variety over an algebraically closed field $K$ of characteristic $p>0$ which is defined over a finite field $\mathbb{F}_{q}$. Let $F=\mathrm{Frob}_{q^2}$ be a Frobenius endomorphism of $G$. Suppose there are points $\alpha,\beta\in G(K)$ such that $F(\alpha)=q^{2}\cdot\alpha$ and $F(\beta)=q\cdot\beta$ (this may happen as we may play with tori and supersingular elliptic curves). Let $g=\alpha+\beta$. Let $d\in\mathbb{Z}_+$ and let $c_1,\dots,c_d\in\mathbb{Z}$. If there is a closed subvariety $X\subseteq G$ such that $X(K)\cap(\mathbb{Z}[F]\cdot g)=\{F^{n_1}(c_1\alpha)+\sum\limits_{i=2}^{d}F^{n_i}(c_{i}\alpha+c_{i-1}\beta)+F^{n_{d+1}}(c_{d}\beta)|\ n_1,\dots,n_{d+1}\in\mathbb{N}\}=\{(\sum\limits_{i=1}^{d}c_{i}q^{2n_i})\cdot\alpha+(\sum\limits_{i=1}^{d}c_{i}q^{n_{i+1}})\cdot\beta|\ n_1,\dots,n_{d+1}\in\mathbb{N}\}$, then $X(K)\cap(\mathbb{Z}\cdot g)$ will contain $\{(\sum\limits_{i=1}^{d}c_{i}q^{2^{i}n})\cdot g|\ n\in\mathbb{N}\}$.
\end{example}

A more precise (but still only heuristic) example is as below.

\begin{example}
Let $p,K$, and $(E,O)$ be same as in Example 3.6(i). Let $G=\mathbb{G}_{m}^{2}\times E^{2}$. Let $C_1,C_2,C_3$ be curves in $G$ given by $C_1(K)=\{(x+1,x,O,O)|\ x\in K\},C_2(K)=\{(1,1,(y+1,\pm\sqrt{(y+1)^{3}+1}),\\(y,\pm\sqrt{y^{3}+1}))|\ y\in K\}$ and $C_3(K)=\{(z+1,z,(z+1,\pm\sqrt{(z+1)^{3}+1}),(z,\pm\sqrt{z^{3}+1}))|\ z\in K\}$ (here we abuse some notations since in fact the sets $C_i(K)$ will be larger as $C_i$ are closed subcurves of $G$). Let $X=C_1+C_2+C_3$ be a closed subvariety of $G$ and let $g=(t+1,t,(t+1,\sqrt{(t+1)^{3}+1}),\\(t,\sqrt{t^{3}+1}))\in G(K)$. Then one can see that $\{(p^n+p^{2n})\cdot g\in X(K)|\ n\in\mathbb{N}\}\subseteq X(K)$ because $(t^{p^{n}}+1,t^{p^n},O,O)\in C_1(K),(1,1,p^{2n}\cdot(t+1,\sqrt{(t+1)^{3}+1}),p^{2n}\cdot(t,\sqrt{t^{3}+1}))\in C_2(K)$ and $(t^{p^{2n}}+1,\\t^{p^{2n}},p^{n}\cdot(t+1,\sqrt{(t+1)^{3}+1}),p^{n}\cdot(t,\sqrt{t^{3}+1}))\in C_3(K)$ for every $n\in\mathbb{N}$.
\end{example}

We believe that this example should already serve as a counterexample of Conjecture 5.2. But as we have mentioned in Remark 3.7, we are unable to completely determine the return set in this example. So we will modify this example and give a rigorous disproof of Conjecture 5.2 by contradiction in next subsection.

\subsection{A rigorous disproof of Conjecture 5.2}
In this subsection, we will prove that Conjecture 5.2 fails for certain translation of isotrivial semiabelian varieties. We fix the data $p,K,(E,O)$ as in Example 3.6(i) through this subsection. Namely, we let our prime $p=5$, let the base field $K$ be $\overline{\mathbb{F}_{p}(t)}$ and let $(E,O)$ be that supersingular elliptic curve as in Example 3.6(i). Our main statement is as follows.

\begin{proposition}
Let $G=\mathbb{G}_{m}^{3}\times(\mathbb{G}_{m}^{2}\times E^{2})$. Let
$$g=\left(\left(t+\alpha,t-\alpha,t\right),\left(t+1,t,(t+1,\sqrt{(t+1)^{3}+1}),(t,\sqrt{t^{3}+1})\right)\right)\in G(K)$$
in which $\alpha\in\mathbb{F}_{p^2}^{\times}$ is a generator of this cyclic group. Let $V_0\subseteq\mathbb{G}_{m}^{3}$ be the closed subvariety given by the equation $x+y=2z+2\alpha^2$. Then there exists a closed subvariety $V\subseteq\mathbb{G}_{m}^{2}\times E^{2}$ such that $\{n\in\mathbb{N}|\ n\cdot g\in(V_0\times V)(K)\}$ is \emph{not} a $p$-normal set in $\mathbb{N}$.
\end{proposition}

Firstly, we need a lemma which is a direct consequence of $\cite[\text{Corollary 9.4}]{Der07}$. Although all numbers in the reference are assumed to be rational, one can verify that this condition is not used in the proof.

\begin{lemma}
Let $d\in\mathbb{N},q\in\mathbb{R}_{>1}$ and $c_1,c_2,e_0,e_1,\dots,e_d\in\mathbb{R}$. Then the set\\
$$\{(n_1,n_2)\in\mathbb{N}^{2}|\ \exists m_1,m_2,\dots,m_d\in\mathbb{N}\text{ s.t. }c_1q^{n_1}+c_2q^{n_2}=e_0+\sum\limits_{i=1}^{d}e_iq^{m_i}\}$$
is a finite union of sets of the 5 forms below:
\begin{enumerate}
\item
a singleton $\{(n_1,n_2)\}$ for some $n_1,n_2\in\mathbb{N}$.
\item
$\{(n+n_1,n_2)|\ n\in\mathbb{N}\}$ for some $n_1,n_2\in\mathbb{N}$.
\item
$\{(n_1,n+n_2)|\ n\in\mathbb{N}\}$ for some $n_1,n_2\in\mathbb{N}$.
\item
$\{(n+n_1,n+n_2)|\ n\in\mathbb{N}\}$ for some $n_1,n_2\in\mathbb{N}$.
\item
$\{(n_1+n_{10},n_2+n_{20})|\ n_1,n_2\in\mathbb{N}\}$ for some $n_{10},n_{20}\in\mathbb{N}$.
\end{enumerate}
\end{lemma}

Next, we introduce a notation which will be used in the proof of Proposition 5.5. We denote $p_0=p^2=25$.

\begin{definition}
Let $q\in\{p_0^{n}|\ n\in\mathbb{Z}_+\},q_1,q_2\in\{p_0^{n}|\ n\in\mathbb{N}\}$ and $c_0\in\mathbb{N},c_1\in\mathbb{Z}_+$. We denote $A(q;q_1,q_2)$ as the set $\{q_1q^{n_1}+q_2q^{n_2}|\ n_1,n_2\in\mathbb{N}\}$ and denote $B(q;c_0,c_1)$ as the set $\{c_0+c_1q^{n}|\ n\in\mathbb{N}\}$. We will obey the convention that all of these coefficients must lie in their ``domain of definition" (i.e. $q\in\{p_0^{n}|\ n\in\mathbb{Z}_+\},q_1,q_2\in\{p_0^{n}|\ n\in\mathbb{N}\}$ and $c_0\in\mathbb{N},c_1\in\mathbb{Z}_+$) when we use this notation.
\end{definition}

Now we can prove Proposition 5.5.

\proof[Proof of Proposition 5.5]
Assume by contradiction that $\{n\in\mathbb{N}|\ n\cdot g\in(V_0\times V)(K)\}$ is a $p$-normal set in $\mathbb{N}$ for every closed subvariety $V\subseteq\mathbb{G}_{m}^{2}\times E^{2}$. Notice that $\{n\in\mathbb{N}|\ n\cdot(t+\alpha,t-\alpha,t)\in V_0(K)\}=\{p_0^{n_1}+p_0^{n_2}|\ n_1,n_2\in\mathbb{N}\}$. Firstly, we prove that up to a finite set, $\{n\in\mathbb{N}|\ n\cdot g\in(V_0\times V)(K)\}$ is a union of finitely many sets of the form $A(q;q_1,q_2)$ along with finitely many sets of the form $B(q;c_0,c_1)$ in which $V\subseteq\mathbb{G}_{m}^{2}\times E^{2}$ is a closed subvariety and the sets $A(q;q_1,q_2),B(q;c_0,c_1)$ are as in Definition 5.7.

~

For a closed subvariety $V\subseteq\mathbb{G}_{m}^{2}\times E^{2}$, we denote $S(V)$ as the set $\{n\in\mathbb{N}|\ n\cdot g\in(V_0\times V)(K)\}$. So by assumption, all of the sets $S(V)$ are $p$-normal sets in $\mathbb{N}$. But since $S(V)\subseteq A(p_0;1,1)$, we firstly conclude that up to a finite set, $S(V)$ is a finite union of sets of the form $S_{q,d,0}(\frac{c_{0}}{q-1};\frac{c_{i}}{q-1})\cap\mathbb{N}$ as in Definition 1.1. Now by considering the ``least common power" $q_0$ of $p_0$ and all the powers of $p$ (i.e. the number $q$) involved here and split $S_{q,d,0}(\frac{c_{0}}{q-1};\frac{c_{i}}{q-1})$ into a finite union of sets of the form $S_{q_0,d,0}(\frac{c_{0}'}{q_0-1};\frac{c_{i}'}{q_0-1})$, we may conclude that there is an element $q_0\in\{p_0^{n}|\ n\in\mathbb{Z}_+\}$ such that up to a finite set, $S(V)$ is a finite union of sets of the form $S_{q_0,d,0}(\frac{c_{0}}{q_0-1};\frac{c_{i}}{q_0-1})\cap\mathbb{N}$. Now using the fact $S(V)\subseteq A(p_0;1,1)$ once more, we know that up to a finite set, $S(V)$ is a finite union of sets of the form $S_{q_0,d,0}(\frac{c_{0}}{q_0-1};\frac{c_{i}}{q_0-1})\cap A(p_0;1,1)$. But we have $A(p_0;1,1)=\bigcup\limits_{a=0}^{M-1}\bigcup\limits_{b=0}^{M-1}A(q_0;p_0^{a},p_0^{b})$ if $q_0=p_0^{M}$ in which $M\in\mathbb{Z}_{+}$. So by Lemma 5.6 we may conclude that up to a finite set, $S(V)$ is a union of finitely many sets of the form $A(q;q_1,q_2)$ along with finitely many sets of the form $B(q;c_0,c_1)$ for every closed subvariety $V\subseteq\mathbb{G}_{m}^{2}\times E^{2}$.

~

Now let $C_1,C_2,C_3$ be locally closed subvarieties of $\mathbb{G}_{m}^{2}\times E^{2}$ such that for any closed point $u\in(\mathbb{G}_{m}^{2}\times E^{2})(K)$, we have
\begin{enumerate}
\item
$u\in C_1(K)$ if and only if $u=(x+1,x,O,O)$ for some $x\in K\backslash\{0,-1\}$,
\item
$u\in C_2(K)$ if and only if $u=(1,1,(y+1,\pm\sqrt{(y+1)^{3}+1}),(y,\pm\sqrt{y^{3}+1}))$ for some $y\in K$, and
\item
$u\in C_3(K)$ if and only if $u=(z+1,z,(z+1,\pm\sqrt{(z+1)^{3}+1}),(z,\pm\sqrt{z^{3}+1}))$ for some $z\in K\backslash\{0,-1\}$.
\end{enumerate}

In fact $C_1$ is a closed subvariety of $\mathbb{G}_{m}^{2}\times E^{2}$ but $C_2,C_3$ are just locally closed subvarieties. Here we use ``$\pm$" to remind that there are two square roots in $K$ while in the definition of $g\in G(K)$, we just arbitrarily fix one. So for example $\pm\sqrt{y^{3}+1}$ should be comprehended as ``an element in $K$ whose square is $y^3+1$".

Let $X$ be the image of the composition map $C_1\times C_2\times C_3\hookrightarrow(\mathbb{G}_{m}^{2}\times E^{2})^{3}\stackrel{m}\rightarrow\mathbb{G}_{m}^{2}\times E^{2}$ where $m$ is the addition map $(u,v,w)\mapsto u+v+w$. Then $X$ is a constructible set in $\mathbb{G}_{m}^{2}\times E^{2}$ satisfying that for any closed point $u\in(\mathbb{G}_{m}^{2}\times E^{2})(K)$, we have $u\in X$ if and only if there exist $x,z\in K\backslash\{0,-1\}$ and $y\in K$ such that $u=((x+1)(z+1),xz,(y+1,\pm\sqrt{(y+1)^{3}+1})+(z+1,\pm\sqrt{(z+1)^{3}+1}),(y,\pm\sqrt{y^{3}+1})+(z,\pm\sqrt{z^{3}+1}))$.

Write $X=\bigcup\limits_{i=1}^{N}(V_i\backslash W_i)$ in which $V_1,\dots,V_N,W_1,\dots,W_N\subseteq\mathbb{G}_{m}^{2}\times E^{2}$ are closed subvarieties satisfying $W_i\subseteq V_i$ for every $1\leq i\leq N$. So $\bigcup\limits_{i=1}^{N}(S(V_i)\backslash S(W_i))=\{n\in\mathbb{N}|\ n\cdot(t+\alpha,t-\alpha,t)\in V_0(K)\}\cap\{n\in\mathbb{N}|\ n\cdot g_0\in X\}$ in which $g_0=(t+1,t,(t+1,\sqrt{(t+1)^{3}+1}),(t,\sqrt{t^{3}+1}))\in(\mathbb{G}_{m}^{2}\times E^{2})(K)$. As a result, we can write $\bigcup\limits_{i=1}^{N}(S(V_i)\backslash S(W_i))=A(p_0;1,1)\cap\{n\in\mathbb{N}|\ \exists x,z\in K\backslash\{0,-1\},\exists y\in K,\ s.t.\\n\cdot g_0=((x+1)(z+1),xz,(y+1,\pm\sqrt{(y+1)^{3}+1})+(z+1,\pm\sqrt{(z+1)^{3}+1}),(y,\pm\sqrt{y^{3}+1})+(z,\pm\sqrt{z^{3}+1}))\}$. Then we know $\{p_0^{n}+p_0^{2n}|\ n\in\mathbb{N}\}\subseteq\bigcup\limits_{i=1}^{N}(S(V_i)\backslash S(W_i))$ by the calculation in Example 5.4. We shall prove that there exists $q_0\in\{p_0^{n}|\ n\in\mathbb{Z}_+\}$ and $q_{10},q_{20}\in\{p_0^{n}|\ n\in\mathbb{N}\}$ such that $\{(q_{10}+q_{20})q_0^{n}|\ n\in\mathbb{N}\}\subseteq\bigcup\limits_{i=1}^{N}(S(V_i)\backslash S(W_i))$.

~

Firstly, fix an element $i_0\in\{1,\dots,N\}$ such that $S(V_{i_0})\backslash S(W_{i_0})$ contains infinitely many elements of the form $p_0^{n}+p_0^{2n}$. Recall that up to a finite set, $S(V)$ is a union of finite sets of the form $A(q;q_1,q_2)$ along with finite sets of the form $B(q;c_0,c_1)$ for every closed subvariety $V\subseteq\mathbb{G}_{m}^{2}\times E^{2}$. So by using the trick of considering the least common power of all the $q$ involved here, we may fix an element $q_0\in\{p_0^{n}|\ n\in\mathbb{Z}_+\}$ and find $E,F\subseteq\mathbb{N}$ such that
\begin{enumerate}
\item
up to a finite set, $S(V_{i_0})$ is equal to $E$ and $S(W_{i_0})$ is equal to $F$, and
\item
both $E$ and $F$ are a union of finite sets of the form $A(q_0;q_1,q_2)$ along with finite sets of the form $B(q_0;c_0,c_1)$.
\end{enumerate}

So in particular the sets $S(V_{i_0})\backslash S(W_{i_0})$ and $E\backslash F$ are equal up to a finite set and hence $E\backslash F$ contains infinitely many elements of the form $p_0^{n}+p_0^{2n}$. But a set of the form $B(q_0;c_0,c_1)$ can only contain finitely many elements of the form $p_0^{n}+p_0^{2n}$. So we conclude that there exists a set $A(q_0;q_{10}',q_{20}')\subseteq E$ such that $A(q_0;q_{10}',q_{20}')\backslash F$ contains infinitely many elements of the form $p_0^{n}+p_0^{2n}$.

Now using Lemma 5.6, we can see that this condition forces $A(q_0;q_{10}',q_{20}')\cap F$ to be a union of a finite set with finitely many sets of the form $B(q_0;c_0,c_1)$. So for a sufficiently large positive integer $C$, we have $\{(q_{10}'q_0^{C}+q_{20}')q_{0}^{n}|\ n\in\mathbb{N}\}\subseteq A(q_0;q_{10}',q_{20}')$ and $\{(q_{10}'q_0^{C}+q_{20}')q_{0}^{n}|\ n\in\mathbb{N}\}\cap F$ is a finite set. So up to a finite set, $\{(q_{10}'q_0^{C}+q_{20}')q_{0}^{n}|\ n\in\mathbb{N}\}$ is contained in $E\backslash F$ and hence it is also contained in $S(V_{i_0})\backslash S(W_{i_0})$ up to a finite set. Therefore, there exists a positive integer $N_0$ such that $\{(q_{10}'q_0^{C+N_0}+q_{20}'q_0^{N_0})q_{0}^{n}|\ n\in\mathbb{N}\}\subseteq S(V_{i_0})\backslash S(W_{i_0})$ and thus we have proved that there exists $q_0\in\{p_0^{n}|\ n\in\mathbb{Z}_+\}$ and $q_{10},q_{20}\in\{p_0^{n}|\ n\in\mathbb{N}\}$ such that $\{(q_{10}+q_{20})q_0^{n}|\ n\in\mathbb{N}\}\subseteq\bigcup\limits_{i=1}^{N}(S(V_i)\backslash S(W_i))$.

~

Now we have $\{(q_{10}+q_{20})q_0^{n}|\ n\in\mathbb{N}\}\subseteq\{n\in\mathbb{N}|\ \exists x,z\in K\backslash\{0,-1\},\exists y\in K,s.t.\ n\cdot g_0=((x+1)(z+1),xz,(y+1,\pm\sqrt{(y+1)^{3}+1})+(z+1,\pm\sqrt{(z+1)^{3}+1}),(y,\pm\sqrt{y^{3}+1})+(z,\pm\sqrt{z^{3}+1}))\}$. This means that for every $n\in\mathbb{N}$, the set of equations
$$
\left\{
\begin{array}{cc}
(t+1)^{(q_{10}+q_{20})q_0^{n}}=(x+1)(z+1) \\
t^{(q_{10}+q_{20})q_0^{n}}=xz \\
(q_{10}+q_{20})q_0^{n}\cdot(t+1,\sqrt{(t+1)^{3}+1})=(y+1,\pm\sqrt{(y+1)^{3}+1})+(z+1,\pm\sqrt{(z+1)^{3}+1}) \\
(q_{10}+q_{20})q_0^{n}\cdot(t,\sqrt{t^{3}+1})=(y,\pm\sqrt{y^{3}+1})+(z,\pm\sqrt{z^{3}+1})
\end{array}
\right.
$$
has a solution $(x_n,y_n,z_n)\in K^{3}$. But since both $q_{10}q_0^{n}$ and $q_{20}q_0^{n}$ are powers of $p$, the first two equations forces $\{x_n,z_n\}=\{t^{q_{10}q_0^{n}},t^{q_{20}q_0^{n}}\}$ and hence we may assume that there exist infinitely many $n\in\mathbb{N}$ such that $x_n=t^{q_{10}q_0^{n}}$ and $z_n=t^{q_{20}q_0^{n}}$ without loss of generality.

Now for every such $n\in\mathbb{N}$, the last two equations can be read as
$$
\left\{
\begin{array}{cc}
((q_{10}+q_{20})q_0^{n}\pm\sqrt{q_{20}q_0^{n}})\cdot(t+1,\sqrt{(t+1)^{3}+1})=(y_n+1,\pm\sqrt{(y_n+1)^{3}+1}) \\
((q_{10}+q_{20})q_0^{n}\pm\sqrt{q_{20}q_0^{n}})\cdot(t,\sqrt{t^{3}+1})=(y_n,\pm\sqrt{y_n^{3}+1})
\end{array}
\right.
$$
for some suitable choice of the sign $\pm$ since we have $\mathrm{Frob}_{p}^{2}=[-p]$ on $E$ (notice that $\sqrt{q_{20}q_0^{n}}$ is a positive integer because both $q_0$ and $q_{20}$ are powers of $p_0=p^2$). But by the explicit formula in $\cite[\text{(III, Ex. 3.7)}]{Sil09}$, one can see that the two signs $\pm$ before $\sqrt{q_{20}q_0^{n}}$ in the two equations must be same. So without loss of generality we assume that there exist infinitely many $n\in\mathbb{N}$ such that
$$
\left\{
\begin{array}{cc}
((q_{10}+q_{20})q_0^{n}+\sqrt{q_{20}q_0^{n}})\cdot(t+1,\sqrt{(t+1)^{3}+1})=(y_n+1,\pm\sqrt{(y_n+1)^{3}+1}) \\
((q_{10}+q_{20})q_0^{n}+\sqrt{q_{20}q_0^{n}})\cdot(t,\sqrt{t^{3}+1})=(y_n,\pm\sqrt{y_n^{3}+1})
\end{array}
\right.
$$

Now let $S\subseteq\mathbb{N}$ be the return set in Example 3.6(i). By Example 3.6, we know that
\begin{enumerate}
\item
$S$ contains infinitely many positive integers of the form $(q_{10}+q_{20})q_0^{n}+\sqrt{q_{20}q_0^{n}}$,
\item
$S\subseteq\{0\}\cup\{p^{k}m|\ k\in\mathbb{N},m\in\mathbb{Z}_{+},m\equiv\pm1\ (\mathrm{mod}\ 2p)\}$, and
\item
$S$ is a union of finitely many arithmetic progressions in $\mathbb{N}$ along with finitely many sets of the form $S_{q,1,0}(\frac{c_0}{q-1};\frac{c_1}{q-1})\cap\mathbb{N}$ where $q$ is a power of $p$ and $c_0,c_1$ are integers satisfying $q-1\mid c_0+c_1$.
\end{enumerate}

But since $p=5>3$, condition (ii) implies that for any arithmetic progression $\{a+dn|\ n\in\mathbb{N}\}\subseteq S$ in which $a\in\mathbb{N}$ and $d\in\mathbb{Z}_+$, we have $v_p(d)>v_p(a+dn)$ for any $n\in\mathbb{N}$ where $v_p(x)$ is the number such that $p^{v_p(x)}\parallel x$. So every arithmetic progression contained in $S$ can only contain finitely many numbers of the form $(q_{10}+q_{20})q_0^{n}+\sqrt{q_{20}q_0^{n}}$. But every set of the form $S_{q,1,0}(\frac{c_0}{q-1};\frac{c_1}{q-1})$ also can only contain finitely many numbers of the form $(q_{10}+q_{20})q_0^{n}+\sqrt{q_{20}q_0^{n}}$. So we get a contradiction.

All in all, we have proved that there exists a closed subvariety $V\subseteq\mathbb{G}_{m}^{2}\times E^{2}$ such that $\{n\in\mathbb{N}|\ n\cdot g\in(V_0\times V)(K)\}$ is not a $p$-normal set in $\mathbb{N}$.
\endproof

One can see that Proposition 5.5 disproves Conjecture 5.2. In other words, in the setting of Theorem 1.5, one cannot require all of the $r$ in the expression of the return set to be 0. Furthermore, Example 5.3 suggests that one cannot require all of the $r$ to be bounded by an absolute constant. We believe that one may prove this statement rigorously by an argument similar to our argument here.

\section{Appendix: an arithmetic version of the $p$-Mordell--Lang problem}

Let $G$ be a semiabelian variety over an algebraically closed field $K$ of positive characteristic. Hrushovski's groundbreaking work $\cite[\text{Theorem}\ 1.1]{Hru96}$ can be regarded as the geometric version of the $p$-Mordell--Lang problem. It describes the structure of a closed subvariety $X\subseteq G$ which admits a dense intersection with a finitely generated group $\Gamma\subseteq G(K)$. The arithmetic version of the $p$-Mordell--Lang problem should describe the form of the intersection set of a closed subvariety $X$ with a finitely generated subgroup $\Gamma$. But it turns out that this intersection set can be quite wild (see $\cite[\text{Section 2}]{GY}$). So one cannot expect the arithmetic version to have a succinct statement. Some works to this end have been done in $\cite[\text{Theorem 1.6, Theorem 1.10}]{GY}$, but it turns out that there is an error in that work. We fix that error and generalize the results in $\cite{GY}$ to arbitrary algebraic groups here. Our main result is as follows.

\begin{theorem}
Let $G$ be an algebraic group over an algebraically closed field $K$ of characteristic $p>0$. Let $X\subseteq G$ be a closed subvariety and let $\Gamma\subseteq G(K)$ be a finitely generated commutative subgroup. Then $X(K)\cap\Gamma$ is a finite union of sets of the form
$$x_0+(\pi|_{\Gamma_0})^{-1}(S)$$
where

$\bullet\ x_0\in\Gamma$,

$\bullet\ G_0\subseteq G$ is an algebraic subgroup which is a semiabelian variety over $K$, and $\Gamma_0=G_0(K)\cap\Gamma$,

$\bullet\ H_0$ is a semiabelian variety over a finite subfield $\mathbb{F}_{q}\subseteq K$, and $F_0$ is the absolute Frobenius endomorphism of $H_0$ corresponding to $\mathbb{F}_{q}$,

$\bullet\ H=H_{0}\times_{\mathbb{F}_{q}}K$, and $F=F_{0}\times_{\mathbb{F}_{q}}K$ is the Frobenius endomorphism of $H$,

$\bullet\ \pi:G_0\rightarrow H$ is a surjective algebraic group homomorphism, and

$\bullet\ S$ is a subset of $\pi(\Gamma_0)$ of the form $\{\alpha_0+\sum\limits_{i=1}^{d}\sum\limits_{j=0}^{r} F^{2^{j}n_{i}}(\alpha_i)|\ n_{1},\dots,n_{d}\in\mathbb{N}\}$ where $d\in\mathbb{Z}_{+},r\in\mathbb{N}$ and $\alpha_0,\alpha_1,\dots,\alpha_{d}\in H(K)$.

Note that since $\Gamma$ is assumed to be commutative, we use $``+"$ for the multiply operation of elements in $\Gamma$.
\end{theorem}

We shall just briefly sketch the proof of Theorem 6.1 since the details are in fact essentially same as the corresponding arguments in $\cite{GY}$ and Section 3 above. Firstly, one can prove Theorem 6.1 in the case of semiabelian varieties in a way just the same as the proof of $\cite[\text{Theorem 1.6, Theorem 1.10}]{GY}$. The only thing to change is that one should use Remark 2.7(ii) above instead of $\cite[\text{Theorem 3.1}]{Ghi08}$ in the proof. This explains where the sets $\{\alpha_0+\sum\limits_{i=1}^{d}\sum\limits_{j=0}^{r} F^{2^{j}n_{i}}(\alpha_i)|\ n_{1},\dots,n_{d}\in\mathbb{N}\}$ in Theorem 6.1 come from. In order to prove the statement for general algebraic groups, one may inherit the strategy in Section 3. The key point is that one can reduce the general case to the case of semiabelian varieties by using $\cite[\text{Theorem 5.6.3(i)}]{Bri17}$ as in the proof of Lemma 3.5 and Theorem 3.1.

\section*{Funding}
This work was supported by the National Natural Science Foundation of China [Grant No. 12271007 to J.X.].

\section*{Acknowledgements}
We thank Dragos Ghioca for comments on the first version of this article and we would like to thank Chengyuan Yang for some useful discussions.

\bibliographystyle{alpha}

\end{spacing}
\end{document}